                      \def\version{30 March 2006}                          %

\documentclass[reqno,11pt]{amsart} 
\usepackage{amsmath} 
\usepackage{amssymb} 
 
\def\emptyset{\varnothing} 

\def\d{\delta} 
 
\def\l{\lambda} 

\def\L{\Lambda}

\newfam\Bbbfam 
\font\tenBbb=msbm10 
\font\sevenBbb=msbm7 
\font\fiveBbb=msbm5 
\textfont\Bbbfam=\tenBbb 
\scriptfont\Bbbfam=\sevenBbb 
\scriptscriptfont\Bbbfam=\fiveBbb

\newcommand{\R}     {\mathbb{R}} 
 
\newcommand{\N}     {\mathbb{N}} 
\renewcommand{\P}   {\mathbb{P}} 
 
\newcommand{\E}     {\mathbb{E}} 
 
 \newcommand{\floor}[1]{\left\lfloor #1 \right\rfloor}

\def\1{{\mathchoice {1\mskip-4mu\mathrm l}      
{1\mskip-4mu\mathrm l} 
{1\mskip-4.5mu\mathrm l} {1\mskip-5mu\mathrm l}}} 
\newcommand{\ssup}[1] {{\scriptscriptstyle{({#1}})}} 
\def\comment#1{} 
\newtheoremstyle{thm}{2ex}{2ex}{\itshape\rmfamily}{} 
{\bfseries\rmfamily}{}{1.7ex}{} 
 
\newtheoremstyle{rem}{1.3ex}{1.3ex}{\rmfamily}{} 
{\itshape\rmfamily}{}{1.5ex}{} 
 
\newenvironment{proofsect}[1] 
{\vskip0.1cm\noindent{\bf #1.}\hskip0.5cm}

\newtheorem{theorem}{Theorem}[section] 
\newtheorem{lemma}[theorem]{Lemma} 
\newtheorem{prop}[theorem] {Proposition} 
\newtheorem{cor}[theorem]  {Corollary}

\theoremstyle{definition}

 \newcommand{\s}{\sigma}

\renewcommand{\section}{\secdef\sct\sect} 
\newcommand{\sct}[2][default]{\refstepcounter{section} 
\vspace{0.8cm} 
\setcounter{equation}{0} 
\centerline{ 
\large\scshape \arabic{section}.\ #1} 
\vspace{0.2cm}} 
\newcommand{\sect}[1]{ 
\vspace{0.8cm} 
\centerline{\large\scshape #1} 
\vspace{0.2cm}} 
 
\renewcommand{\subsection}{\secdef \subsct\sbsect} 
\newcommand{\subsct}[2][default]{\refstepcounter{subsection} 
\nopagebreak 
\vspace{0.5\baselineskip} 
{\flushleft\bf \arabic{section}.\arabic{subsection}~\bf #1  } 
\nopagebreak} 
\newcommand{\sbsect}[1]{\vspace{0.1cm}\noindent 
{\bf #1}\vspace{0.1cm}} 
\newcommand{\heap}[2]{\genfrac{}{}{0pt}{}{#1}{#2}} 
 
\renewcommand{\subsubsection}{%
\secdef \subsubsect\sbsbsect} 
\newcommand{\subsubsect}[2][default]{%
\refstepcounter{subsubsection} 
\nopagebreak 
\vspace{0.1\baselineskip} 
\nopagebreak 
{\flushleft 
\sffamily\slshape 
\arabic{section}.\arabic{subsection}.\arabic{subsubsection} 
\ %
\sffamily #1\/.}\ } 
\newcommand{\sbsbsect}[1]{\vspace{0.1cm}\noindent 
{\bf #1}\ } 
 
\renewcommand{\d}{{\rm d}} 
 
\newcommand{\eps}{\varepsilon} 
\newcommand{\Leb}{{\rm Leb}} 
\newcommand{\Sym}{\mathfrak{S}}

\newcommand{\supp}{{\operatorname {supp}}} 
\newcommand{\dist}{{\operatorname {dist}}} 
\newcommand{\diam}{{\operatorname {diam}}}

\newcommand{\Ccal}   {{\mathcal C }}

\newcommand{\Hcal}   {{\mathcal H }}

\newcommand{\Lcal}   {{\mathcal L }} 
\newcommand{\Mcal}   {{\mathcal M }}

\newcommand{\Ucal}   {{\mathcal U }}

\newcommand{\Xcal}   {{\mathcal X }}

\newcommand{\m} {{\mathfrak m}}
 
\begin{document} 
 
\title[Large deviations for symmetrised Brownian bridges]{\large 
Large deviations for many Brownian bridges\\ \medskip with symmetrised initial-terminal condition} 
 
\author[Stefan Adams and Wolfgang 
        K{\"o}nig]{} 
\maketitle

\thispagestyle{empty} 
\vspace{0.2cm} 
 
\centerline {\sc By Stefan Adams\footnote{Dublin Institute for Advanced Studies, School of Theoretical Physics, 10, Burlington Road, Dublin 4, Max-Planck Institute for Mathematics in the Sciences, Inselstra{\ss}e 22-26, D-04103 Leipzig, Germany, {\tt adams@mis.mpg.de}} \/ and  Wolfgang
K{\"o}nig\footnote{Mathematisches Institut, Universit\"at Leipzig, Augustusplatz 10/11, D-04109 Leipzig, Germany, {\tt koenig@math.uni-leipzig.de}}}
\vspace{0.4cm}
\renewcommand{\thefootnote}{}
 \footnote{
Partially supported by DFG grant AD 194/1-1 and by the DFG-Forschergruppe 718 \lq Analysis and stochastics in complex physical systems\rq}

\centerline{\small(\version)} 
\vspace{.5cm} 
 

\bigskip

\begin{quote} 
{\small {\bf Abstract:}} Consider a large system of $N$ Brownian motions in $\R^d$ with some non-degenerate initial measure on some fixed time interval $[0,\beta]$ with symmetrised initial-terminal condition. That is, for any $i$, the terminal location of the $i$-th motion is affixed to the initial point of the $\s(i)$-th motion, where $\s$ is a uniformly distributed random permutation of $1,\dots,N$. Such systems play an important role in quantum physics in the description of Boson systems at positive temperature $1/\beta$.

In this paper, we describe the large-$N$ behaviour of the empirical path measure (the mean of the Dirac measures in the $N$ paths) and of the mean of the normalised occupation measures of the $N$ motions in terms of large deviations principles. The rate functions are given as variational formulas involving certain entropies and Fenchel-Legendre transforms. Consequences are drawn for asymptotic independence statements and laws of large numbers.

In the special case related to quantum physics, our rate function for the occupation measures turns out to be equal to the well-known Donsker-Varadhan rate function for the occupation measures of one motion in the limit of diverging time. This enables us to prove a simple formula for the large-$N$ asymptotic of the symmetrised trace of ${\rm e}^{-\beta \Hcal_N}$, where $\Hcal_N$ is an $N$-particle Hamilton operator in a trap.
\end{quote}

\vfill


 
\eject 
 
\setcounter{section}{0} 
\section{Introduction and main results}\label{Intro}

\subsection{Introduction.}

\noindent In this article, we study the large-$N$ behaviour of a system of $N$ symmetrised Brownian motions in $\R^d$ on a fixed time interval $[0,\beta]$, i.e., the behaviour of the system under the measure
\begin{equation}\label{Psymdef}
\P_{\m,N}^{\ssup{{\rm sym}}}=\frac{1}{N!}\sum_{\s\in\Sym_N}\;\int_{\R^d}\cdots \int_{\R^d} \m(\d x_1)\cdots\m(\d x_N)\bigotimes_{i=1}^N\P^{\beta}_{x_i,x_{\s(i)}}. 
\end{equation}
Here $\Sym_N$ is the set of all permutations of $1,\dots,N$, $\P_{x,y}^{\beta}$ the normalised Brownian bridge measure on the time interval $[0,\beta]$ with initial point $ x \in\R^d$ and terminal point $ y \in\R^d$ (also see \eqref{nnBBM} below), and $\m$ is the initial probability distribution on $\R^d$.  Hence, the terminal location of the $i$-th motion is affixed to the initial location of the $\sigma(i)$-th motion, where $\sigma$ is a uniformly distributed random permutation. Any of the $N$ paths is a Brownian motion with initial distribution $\m$, but with a peculiar terminal distribution at time $\beta$. We can conceive $\P_{\m,N}^{\ssup{{\rm sym}}}$ as a two-step random mechanism: First we pick a uniform random permutation $\sigma$, then we pick $N$ Brownian motions with initial distribution $\m$, and the $i$-th motion is conditioned to terminate at the initial point of the $\sigma(i)$-th motion, for any $i$.

One main motivation to consider this model stems from quantum physics, where one is interested in the description of the canonical ensemble of large Boson systems at positive temperature, see Section~\ref{sec-Quantum} below. Beside the application in physics, the problem is also appealing from a mathematical point of view, since the combinatorics of a random permutation is combined with independent, but not identically distributed, objects.

We consider the distribution of the tuple of $N$ random paths $B^{\ssup 1},\dots,B^{\ssup N}\colon[0,\beta]\to\R^d$ under $\P_{\m,N}^{\ssup{{\rm sym}}}$. We are interested in the large-$N$ behaviour of the empirical path measure
\begin{equation}\label{LNdef}
L_N=\frac{1}{N}\sum_{i=1}^N\delta_{B^{\ssup{i}}}\in\Mcal_1(\Ccal),
\end{equation}
which is a random probability measure on the set $ \Ccal$ of continuous paths $[0,\beta]\to\R^d$. More precisely, we derive a large deviations principle for the distributions of $L_N$ under $\P_{\m,N}^{\ssup{{\rm sym}}}$ as $N\to\infty$ (Theorem~\ref{main}). (In Section~\ref{sec-LDP} below we recall the notion of a large deviations principle.) We also obtain a large deviations principle for the means of the normalised occupation measures,
\begin{equation}\label{YNdef}
Y_N= \frac{1}{N}\sum_{i=1}^N\frac{1}{\beta}\int_0^{\beta}\d s\,\delta_{B^{\ssup{i}}_s}\in\Mcal_1(\R^d),
\end{equation}
(Theorem~\ref{main2}).
Our large-deviation rate functions for the two principles are explicit in terms of variational problems involving an entropy term (describing the large deviations of the permutations) and a certain Legendre transform (describing the large deviations of $L_N$ and $Y_N$, respectively, for a fixed permutation). We draw a number of corollaries about variants of the principles, laws of large numbers and asymptotic independence.
 
Our results are most beautiful and most striking for the important special case that $\m$ is the Lebesgue measure on a bounded box and that $\P_{x,y}^\beta$ is replaced by the canonical, non-normalised, Brownian bridge measure, $\mu_{x,y}^\beta$ (see \eqref{nnBBM}). In this case, the rate function for the means $Y_N$ turns out to be equal to $ \beta $ times the well-known Donsker-Varadhan rate function, which is explicitly given as the energy of the square root of the density of the measure considered (Theorem~\ref{Jident}). This function is well-known as the rate function for the normalised occupation measure for just one Brownian motion (or bridge) in the limit of diverging time. We give an interpretation of this remarkable coincidence in terms of the well-known cycle structure of the permutations in \eqref{Psymdef}. However, let us remark that our proofs do not respect this structure at all and therefore give no rigorous insight into that.

The mentioned identification of the rate function for the $Y_N$'s as the Donsker-Varadhan rate function has interesting consequences for the asymptotic description of the canonical ensemble of a system of $N$ noninteracting Bosons at positive temperature $1/\beta\in(0,\infty)$. In fact, we obtain in Section~\ref{sec-Quantum} a remarkably simple formula for the large-$N$ asymptotic of the symmetrised trace of the $N$-particle Hamilton operator in a fixed box (Theorem~\ref{LDPW}). We consider this as a first step towards a rigorous understanding of large Boson systems at positive temperature. Future work will be devoted to the mutually interacting case. Interacting Brownian motions in trap potentials have been so far analysed without symmetrisation, in particular, finite systems for vanishing temperature in \cite{ABK05a} and large systems of interacting motions for fixed positive temperature in \cite{ABK05b}.

Let us make some remarks on related literature. In \cite{S31} Schr\"odinger raised the question of the most probable behaviour of a large system of diffusion particles in thermal equilibrium. F\"ollmer \cite{F88} gave a mathematical formulation of theses ideas in terms of large deviations. He applied Sanov's theorem to obtain a large deviations principle for $ L_N $ when $ B^{\ssup{1}},B^{\ssup{2}},\ldots $ are i.i.d.~Brownian motions with initial distribution $\m$ and no condition at time $\beta$. The rate function is the relative entropy with respect to $ \int_{\R^d}\m(\d x)\,\P_{x}\circ B^{-1} $, where the motions start in $x$ under $\P_x$. Then Schr\"odinger's question amounts to identifying the minimiser of that rate function under given fixed independent initial and terminal distributions. Interestingly, it turns out that the unique minimiser is of the form $\int_{\R^d}\int_{\R^d}\d x\d y\,f(x) g(y)\,\P_{x,y}^\beta\circ B^{-1}$, i.e., a Brownian bridge with independent initial and terminal distributions. The probability densities $f$ and $g$ are characterised by a pair of dual variational equations, originally appearing in \cite{S31} for the special case that both given initial and terminal measures are the Lebesgue measure. The monograph \cite{N93} systematically studies such dual equations and their connections to the dual time-dependent Schr\"odinger equations and to  Schr\"odinger processes, i.e., processes of the form $\int_{\R^d}\int_{\R^d}q(\d x,\d y)\,\P_{x,y}^\beta$ that are additionally Markov. \cite{FG97} obtained conditions from minimising the entropy  to derive the Markov property of such processes; the absolute continuity and the product structure of $q$ turn out to be crucial.

An important work combining combinatorics and large deviations for symmetrised measures is \cite{Toth90}. T\'oth \cite{Toth90} considers $ N $ continuous-time simple random walks on a complete graph with $ \floor{\rho N} $ vertices, where $ \rho\in (0,1) $ is fixed. He looks at the symmetrised distribution as in \eqref{Psymdef} and adds an exclusion constraint: there is no collision of any two particles during the time interval $ [0,\beta] $. The combinatorial structure of this model enabled him to express the free energy in terms of a cleverly chosen Markov process on $ \N_0 $. Using Freidlin-Wentzell theory, he derives an explicit formula for the large-$ N $ asymptotic of the free energy; in particular he obtains a phase-transition, called {\it Bose-Einstein-condensation\/}, for large $ \beta $ and sufficiently large $ \rho $.

Our proof is partly inspired by the method developed in \cite{KM02}. The problem there is the evaluation of the large-$k$ asymptotic of the $k$-th moments of the intersection local time in $U$ of $p\in\N\setminus\{1\}$ Brownian motions running in an open subset $O$ of $\R^2$ or $\R^3$. This moment is known to be equal to
$$
\int_U \d x_1\cdots\int_U \d x_k \, \Big(\sum_{\sigma\in \Sym_k} \prod_{i=1}^k G\big(x_{\sigma(i-1)}, x_{\sigma(i)} \big)\Big)^p.
$$
Here $U$ is a compact subset of $O$, and $G$ is the Green's function of one of the Brownian motions in $O$. Even though the motivation for studying this problem is quite different, similar techniques prove useful for the study of that problem and the one of the present paper.

The effect of mixing random variables using a random uniformly distributed permutation to large deviation principles has been studied both in \cite{DZ92} and \cite{T02}, which were motivated from asymptotic questions about exchangeable vectors of random variables.  \cite{DZ92} studies large deviations for the empirical measures $\frac 1N\sum_{i=1}^N\delta_{Y_i}$, where $Y_1,\dots,Y_N$ have distribution $\int_\Theta \mu(\d \theta)\, P_N^{\ssup \theta}$ for some distribution $\mu$ on some compact space $\Theta$, and the empirical measures are assumed to satisfy a large deviation principle under $P_N^{\ssup \theta}$ for each $\theta$. In \cite{T02}, a similar problem is studied: given a sequence of random vectors $(Y_1^{\ssup N},\dots,Y_N^{\ssup N})$ such that the empirical measures $\frac 1N\sum_{i=1}^N\delta_{Y_i^{\ssup N}}$ satisfy a large deviation principle, another principle is established for the process of empirical measures $\frac 1N\sum_{i=1}^{\lfloor tN\rfloor}\delta_{X_{i}^{\ssup N}}$, where 
$$
\big(X_1^{\ssup N},\dots,X_N^{\ssup N}\big)=\frac1{N!}\sum_{\s\in\Sym_N}\big(Y_{\s(1)}^{\ssup N},\dots,Y_{\s(N)}^{\ssup N}\big).
$$
In both works, the large-deviation rate function is expressed in terms of entropy terms, like in our main result. However, a substantial difference is that the symmetrisation mechanism in \eqref{Psymdef} is adequately described only by the pairs $(i,{\s(i)})$ for $i=1,\dots,N$, instead of just the sequence of the $\s(i)$'s.

Let us also mention that our problem may also be seen as a particular two-level large deviations result, which has something in common with general multilevel large deviations as studied in \cite{DG94}. There a large deviation principle is established for $\frac1 M\sum_{i=1}^M\delta_{X_i^{\ssup N}}$, as $M,N\to\infty$, under the assumption that, for any $i$, the sequence $(X_i^{\ssup N})_{N\in\N}$ satisfies a principle.

The structure of the remainder of this paper is the following. In Section~\ref{sec-LDresults} we describe our main results. A couple of remarks and consequences are in Section~\ref{sec-remarks}, and in Section~\ref{sec-spcase} we consider an important particular case, where we identify the rate function in simple terms. Some remarks on relations to quantum physics are in Section~\ref{sec-Quantum}. In Section~\ref{sec-JLambdaident} we prove some facts about the rate functions, and Section~\ref{sec-mainproof} contains the proof of the large deviations principles. Finally, the Appendix, Section~\ref{sec-LDP}, recalls some notions and facts about large deviations theory.

\subsection{Large deviations for $\boldsymbol{\P^{\ssup{\rm sym}}_{\m,N}} $.} \label{sec-LDresults}

\noindent We are going to formulate our first main result: large deviations principles for the distributions of $L_N$ and $ Y_N $ under $\P_{\m,N}^{\ssup{{\rm sym}}}$, see \eqref{LNdef}, \eqref{YNdef} and \eqref{Psymdef}. (We refer to Section~\ref{sec-LDP} for the notion of a large deviation principle and related notation.)  Throughout the paper, we fix $ \beta >0 $.  Let $ \Ccal=\Ccal([0,\beta];\R^d) $ be the set of continuous functions $ [0,\beta]\to \R^d $. We equip $\Ccal$ with the topology of uniform convergence and with the corresponding Borel $\sigma$-field. We consider $N$ random variables, $ B^{\ssup{1}},\ldots, B^{\ssup{N}}$, taking values in $\Ccal$. For the reader's convenience, we repeat the definition of a Brownian bridge measure; see the Appendix in \cite{Sz98}. We decided to work with Brownian motions having generator $\Delta$ instead of $\frac 12\Delta$. We write $\P_x$ for the probability measure under which $B=B^{\ssup{1}}$ starts from $x\in\R^d$. The canonical (non-normalised) Brownian bridge measure on the time interval $[0,\beta]$ with initial site $x\in\R^d$ and terminal site $y\in \R^d$ is defined as
\begin{equation}\label{nnBBM}
\mu_{x,y}^\beta(A)=\frac {\P_x(B\in A;B_\beta\in\d y)}{\d y},\qquad A\subset \Ccal\mbox{ measurable.}
\end{equation}
Then $\mu_{x,y}^\beta$ is a regular Borel measure on $\Ccal$ with total mass equal to the Gaussian density,
\begin{equation}\label{Gaussian}
\mu_{x,y}^\beta(\Ccal)=p_\beta(x,y)=\frac {\P_x(B_\beta\in\d y)}{\d y}=(4\pi\beta)^{-d/2}{\rm e}^{-\frac 1{4\beta}|x-y|^2}.
\end{equation}
The normalised Brownian bridge measure is defined as $\P_{x,y}^\beta=\mu_{x,y}^\beta/p_\beta(x,y)$, which is a probability measure on $\Ccal$.

Now we introduce the rate functions.
By $\Mcal_1(X)$  we denote the set of Borel probability measures on a topological space $X$. With 
\begin{equation}\label{Entropydef}
H(q | \widetilde{q})=\int_{X}q(\d x)\log\frac{q(\d x)}{\widetilde{q}(\d x)}
\end{equation} 
we denote the relative entropy of $q\in \Mcal_1(X)$ with respect to $ \widetilde q\in\Mcal_1(X)$. We will often use this notation for $X=\R^d\times\R^d$ in the sequel, but also for other spaces $X$. Let  $ \Mcal_1^{\ssup{\rm s}}(\R^d\times\R^d)$ be the set of shift-invariant probability measures $q$ on $ \R^d\times\R^d $, i.e., measures whose first and second marginals coincide and are both denoted by $\overline q$. Note that $q\mapsto  H(q|\overline{q}\otimes\m)$ is strictly convex. We write $\langle \Phi ,\mu\rangle$ for integrals $\int_\Ccal \Phi(\omega)\,\mu(\d \omega)$ for suitable functions $\Phi$ on $\Ccal$. Define the functional  $I^{\ssup{{\rm sym}}}_{\m}$ on $ \Mcal_1(\Ccal) $ by
\begin{equation}\label{Isymdef}
I^{\ssup{{\rm sym}}}_{\m}(\mu)=\inf_{q\in\Mcal_1^{\ssup{\rm s}}(\R^d\times\R^d)}\Bigl\{H(q|\overline{q}\otimes\m)+ I^{\ssup{q}}(\mu)\Bigr\},\qquad \mu\in\Mcal_1(\Ccal),
\end{equation} 
where
\begin{equation}\label{Jqdef}
I^{\ssup{q}}(\mu)=\sup_{\Phi \in\Ccal_{\rm b}(\Ccal)}\Bigl\{\langle \Phi ,\mu\rangle -\int_{\R^d}\int_{\R^d} q(\d x, \d y)\log\E_{x,y}^{\beta}\big[{\rm e}^{\Phi(B)}\big]\Bigr\},\qquad\mu\in\Mcal_1(\Ccal).
\end{equation} 
Hence, $I ^{\ssup{q}} $ is a Legendre-Fenchel transform, but {\it not\/} the one of a logarithmic moment generating function of any random variable. In particular, $I^{\ssup{q}}$, and therefore also $I^{\ssup{{\rm sym}}}_{\m}$, are nonnegative, and $I^{\ssup{q}}$ is convex as a supremum of linear functions. There seems to be no way to represent $I^{\ssup{q}}(\mu)$ as the relative entropy of $\mu$ with respect to any measure.

By $\pi_s\colon \Ccal\to \R^d$ we denote the projection $\pi_s(\omega)=\omega_s$. The marginal measure of $\mu\in\Mcal_1(\Ccal)$ is denoted by $\mu_s=\mu\circ\pi_s^{-1}\in\Mcal_1(\R^d)$; analogously we write $\mu_{0,\beta}=\mu\circ(\pi_0,\pi_\beta)^{-1 }\in\Mcal_1(\R^d\times\R^d)$ for the joint distribution of the initial and the terminal point of a random process with distribution $\mu$. It is easy to see that $q=\mu_{0,\beta}$ if  $I ^{\ssup{q}} (\mu)<\infty$. Indeed, in \eqref{Jqdef} relax the supremum over all $\Phi\in\Ccal_{\rm b}(\Ccal)$ to all functions of the form $\omega\mapsto f(\omega_0,\omega_\beta)$ with $f\in\Ccal_{\rm b}(\R^d)$. This gives that
$$
\begin{aligned}
\infty>I ^{\ssup{q}} (\mu)&\geq\sup_{f\in\Ccal_{\rm b}(\R^d)}\Big(\big\langle\mu_{0,\beta},f\big\rangle-\big\langle q,\log \E_{\pi_0,\pi_\beta}^{\beta}\big[{\rm e}^{f(B_0,B_\beta)}\big]\Big)
&=\sup_{f\in\Ccal_{\rm b}(\R^d)}\big\langle\mu_{0,\beta}-q,f\big\rangle,
\end{aligned}
$$
and this implies that $\mu_{0,\beta}=q$. In particular, the infimum in \eqref{Isymdef} is uniquely attained at this $q$, i.e.,
\begin{equation}\label{Isymident}
I^{\ssup{{\rm sym}}}_{\m}(\mu)=\begin{cases}H(\mu_{0,\beta}|\mu_0\otimes\m\big)+\sup\limits_{\Phi \in\Ccal_{\rm b}(\Ccal)}\Big\langle\mu, \Phi  -\log\E_{\pi_0,\pi_\beta}^{\beta}\big[{\rm e}^{\Phi(B)}\big]\Big\rangle&\mbox{if }\mu_0=\mu_\beta,\\
+\infty&\mbox{otherwise.}
                               \end{cases}
\end{equation}
In particular, $I^{\ssup{{\rm sym}}}_{\m}$ is convex.

Then our main result reads as follows.

\begin{theorem}[Large deviations for $L_N$]\label{main}
Fix $ \beta \in (0,\infty) $ and assume that the initial distribution $ \m\in\Mcal_1(\R^d) $ has compact support. Then, as $ N\to\infty$, under the symmetrised measure $ \P_{\m,N}^{\ssup{{\rm sym}}} $, the empirical path measures $ L_N $ satisfy a large deviations principle on $ \Mcal_1(\Ccal) $ with speed $ N $ and rate function $ I^{\ssup{{\rm sym}}}_{\m}$.
\end{theorem}

To be more explicit, the stated large deviations principle says that
$$
\lim_{N\to\infty}\frac 1N\log \P_{\m,N}^{\ssup{{\rm sym}}}\big(L_N\in\,\cdot\big)=-\inf_{\mu\in\,\cdot}I^{\ssup{{\rm sym}}}_{\m}(\mu),
$$
in the weak sense, i.e., there is a lower bound for open subsets of $\Mcal_1(\Ccal)$ and an upper bound for closed ones. The proof of Theorem~\ref{main} is in Section~\ref{sec-mainproof}. In Section~\ref{sec-outline} we give an outline of its main idea. The assumption that the initial distribution $ \m\in\Mcal_1(\R^d) $ has compact support is necessary only in the proof of the lower bound. However, we prove the upper bound and the exponential tightness without using this assumption. Our proof does not rely on the Markov property of the Brownian bridge processes; only some continuity is required, see in particular Lemma~\ref{estimationsetendpoints} below.

We also have the analogous result for the mean of the occupation measures, $Y_N$, defined in \eqref{YNdef}.
For formulating this, we define the functional $J^{\ssup{\rm sym}}_{\m}$ on $ \Mcal_1(\R^d) $ by
\begin{equation}\label{Jsymdef}
J^{\ssup{\rm sym}}_{\m}(p)=\inf_{q\in\Mcal_1^{\ssup{\rm s}}(\R^d\times\R^d)}\Bigl\{H(q|\overline{q}\otimes\m)+ J^{\ssup{q}}(p)\Bigr\},\qquad p\in\Mcal_1(\R^d),
\end{equation}
where
\begin{equation}\label{Jsymqdef}
J^{\ssup{q}}(p)=\sup_{f\in\Ccal_{\rm b}(\R^d)}\Bigl\{\beta\langle f,p\rangle -\int_{\R^d}\int_{\R^d} q(\d x,\d y)\,\log\E_{x,y}^{\beta}\Bigl[{\rm e}^{\int_0^{\beta}f(B_s)\,\d s}\Bigr]\Bigr\}.
\end{equation}

\begin{theorem}[Large deviations for $Y_N$]\label{main2}
Fix $ \beta \in (0,\infty) $ and assume that the initial distribution $ \m\in\Mcal_1(\R^d) $ has compact support. Then, as $ N\to\infty $, under the symmetrised measure $ \P_{\m,N}^{\ssup{{\rm sym}}} $, the mean of occupation measures, $ Y_N $, satisfy a large deviations principle on $ \Mcal_1(\R^d) $ with speed $ N $  and rate function $ J^{\ssup{\rm sym}}_{\m} $.
\end{theorem}

The proof of Theorem~\ref{main2} is in Section~\ref{sec-proofYN}. Via the contraction principle \cite[Th.~4.2.1]{DZ98}, the large deviations principle for $Y_N$ in Theorem~\ref{main2} is a consequence of the one for $L_N$ in Theorem~\ref{main}. Indeed, consider $ \Psi\colon\Mcal_1(\Ccal)\to\Mcal_1(\R^d) $ defined by $ \Psi(\mu)=\frac{1}{\beta}\int_0^{\beta}\d s\, \mu\circ\pi_s^{-1} $, where we recall that $ \pi_s(\omega)=\omega(s) $ is the projection. Then $ \Psi $ is continuous and bounded, and $ Y_N=\Psi(L_N) $. Hence, the contraction principle immediately yields the large deviation principle for $ Y_N $ with rate function given by 
\begin{equation}\label{ratefctJtilde}
\widetilde J^{\ssup{\rm sym}}_\m(p)=\inf_{q\in\Mcal_1^{\ssup{\rm s}}(\R^d\times\R^d)}\big\{H(q|\overline{q}\otimes\m)+ \widetilde J^{\ssup{q}}(p)\big\},
\end{equation}
where
\begin{equation}\label{ratefunctioncontract}
\widetilde J^{\ssup{q}}(p)=\inf\big\{I^{\ssup{q}}(\mu)\colon\mu\in\Mcal_1(\Ccal), \Psi(\mu)=p\big\},\qquad p\in\Mcal_1(\R^d).
\end{equation}
In Section~\ref{sec-proofYN} below we will show that $ \widetilde J^{\ssup{q}}= J^{\ssup{q}} $ for any $q\in\Mcal_1^{\ssup{\rm s}}(\R^d\times\R^d)$, which implies that $\widetilde J^{\ssup{\rm sym}}_\m= J^{\ssup{\rm sym}}_\m$ and finishes the proof of Theorem~\ref{main2}. The proof of $ \widetilde J^{\ssup{q}}\geq J^{\ssup{q}} $ is simple and straightforward. However, a direct analytical proof of the complementary inequality,  $ \widetilde J^{\ssup{q}}\leq J^{\ssup{q}} $, seems rather difficult; we prove this indirectly by showing that $\widetilde J^{\ssup{q}}$ and $J^{\ssup{q}}$ govern the same large deviations principle.

Let us give a brief informal interpretation of the shape of the rate functions in \eqref{Isymdef} and \eqref{Jsymdef}. As we have remarked earlier, the symmetrised measure $ \P_{\m,N}^{\ssup{{\rm sym}}} $ arises from a two-step probability mechanism. This is reflected in the representation of the rate function $ I^{\ssup{{\rm sym}}}_{\m}$ in \eqref{Isymdef}: in a peculiar way (which we roughly describe in Section~\ref{sec-outline}), the entropy term $H(q|\overline{q}\otimes\m)$ describes the large deviations of the uniformly distributed random permutation $\sigma$, together with the integration over $\m^{\otimes N}$. The measure $q$ governs a particular distribution of $N$ independent, but not identically distributed, Brownian bridges. Under this distribution, $L_N$ satisfies a large deviations principle with rate function $I^{\ssup{q}}$, which also can be guessed from the G\"artner-Ellis theorem \cite[Th.~4.5.20]{DZ98}. The presence of a two-step mechanism makes impossible to apply this theorem directly to $ \P_{\m,N}^{\ssup{{\rm sym}}} $. 

Let us contrast this to the case of i.i.d.~Brownian bridges $B^{\ssup 1},\dots,B^{\ssup N}$ with starting distribution $ \m $, i.e., we replace $ \P_{\m,N}^{\ssup{{\rm sym}}} $  by $ (\int\m(\d x)\, \P_{x,x}^\beta)^{\otimes N}$. Here the empirical path measure $ L_N $ satisfies a large deviations principle with rate function 
$$
I_{\m}(\mu)=\sup_{\Phi \in\Ccal_{\rm b}(\Ccal)}\Bigl\{\langle \Phi ,\mu\rangle-\log\int_{\R^d}\m(\d x)\, \E_{x,x}^\beta\bigl[{\rm e}^{\Phi(B)}\bigr]\Bigr\},
$$
as follows from an application of Cram\'er's theorem \cite[Theorem~6.1.3]{DZ98}.  Note that $I_{\m}(\mu)$ is the relative entropy of $\mu$ with respect to $\int\m(\d x)\, \P_{x,x}^\beta\circ B^{-1}$. Although there is apparently no reason to expect a direct comparison between the distributions of $L_N$ under  $ \P_{\m,N}^{\ssup{{\rm sym}}} $ and under $ (\int\m(\d x)\, \P_{x,x}^\beta)^{\otimes N}$, the rate functions admit a simple relation: it is easy to see that $ I^{\ssup{q}}\ge I_{\m} $ for the measure $ q(\d x, \d y)=\m(\d x)\delta_x(\d y)\in \Mcal_1^{\ssup{\rm s}}(\R^d\times\R^d)$, since
\begin{eqnarray}
\begin{aligned}
-\int_{\R^d}\int_{\R^d} q(\d x,\d y)\log\E^{\beta}_{x,y}\big[{\rm e}^{\Phi(B)}\big]&\ge -\log\int_{\R^d}\m(\d x)\, \E_{x,x}^\beta\bigl[{\rm e}^{\Phi(B)}\bigr]. 
\end{aligned}
\end{eqnarray} 
In particular, $ I^{\ssup{\rm sym}}_{\m}\ge I_{\m} $.

All the preceding remarks apply to the mean of the occupation measures, $Y_N$, in place of $L_N$.

\subsection{Extensions and remarks.}\label{sec-remarks}

\noindent Let us extend Theorems~\ref{main} and \ref{main2} to some larger class of measures $\P_{\m,N}^{\ssup{{\rm sym}}}$. Obviously, both theorems remain true if the total mass of the initial measure $\m$ is not necessarily equal to one, but positive and finite. (Here we adapt the notion of a large deviations principle accordingly, which is easily done by dropping the requirement that the infimum of the rate function be equal to zero; see Section~\ref{sec-LDP}.) A bit deeper lies the fact that the Brownian bridge measure does not have to be normalised in order that the results hold:

\begin{prop}\label{cormain}
Fix $ \beta \in (0,\infty) $ and assume that $ \m $ is a positive finite measure on $\R^d$ with compact support. Fix some continuous function $g\colon\R^d\times \R^d\to(0,\infty)$ and replace $\P_{x,y}^\beta$ by $g(x,y)\P_{x,y}^\beta$ in the definition \eqref{Psymdef} of $\P_{\m,N}^{\ssup{{\rm sym}}}$. Then 
\begin{enumerate}
\item[(i)] Theorem~\ref{main} remains true. The corresponding rate function is $\mu\mapsto I^{\ssup{\rm sym}}_{\m}(\mu)-\langle \mu_{0,\beta},\log g\rangle$.

\item[(ii)] Theorem~\ref{main2} remains true. The corresponding rate function is
\begin{equation}\label{Jsymgdef}
J^{\ssup{\rm sym}}_{\m,g}(p)=\inf\limits_{q\in\Mcal_1^{\ssup{{\rm s}}}(\R^d\times\R^d)}\Big\{H(q|\overline{q}\otimes\m)+J^{\ssup{q}}(p)-\langle q,\log g\rangle\Big\},\qquad p\in\Mcal_1(\R^d).
\end{equation}
\end{enumerate}
\end{prop}

\begin{proofsect}{Proof} (i) Note that, for any $\s\in\Sym_N$ and any $x_1,\dots,x_N\in\supp(\m)$, with probability one with respect to $\otimes_{i=1}^N\P_{x_i,x_{\s(i)}}^{\beta}$,
$$
\prod_{i=1}^N g(x_i,x_{\s(i)})={\rm e}^{\sum_{i=1}^N\log g(B_0^{\ssup i},B_\beta^{\ssup i})}
={\rm e}^{N\langle L_N,\Phi_g\rangle},
$$
where $\Phi_g(\omega)=\log g(\omega_0,\omega_\beta)$. Since  $\Phi_g$ is bounded and continuous, the principle follows from \cite[Th.~III.17]{dH00}, and the rate function is identified as $\mu\mapsto I^{\ssup{\rm sym}}_{\m}(\mu)-\langle\mu,\Phi_g\rangle=I^{\ssup{\rm sym}}_{\m}(\mu)-\langle \mu_{0,\beta},\log g\rangle$.

\noindent (ii) In the same way as Theorem~\ref{main2} is deduced from Theorem~\ref{main} via the contraction principle \cite[Th.~4.2.1]{DZ98}, (ii) is derived from (i). Indeed, using the principle in (i), the contraction principle implies the desired principle with rate function
$$
\begin{aligned}
p\mapsto&\inf_{\mu\in\Mcal_1(\Ccal)\colon\Psi(\mu)=p}\big\{I^{\ssup{\rm sym}}_{\m}(\mu)-\langle \mu_{0,\beta},\log g\rangle\big\}\\
&=\inf_{q\in\Mcal_1^{\ssup{\rm s}}(\R^d\times\R^d)}\Big\{H(q|\overline{q}\otimes\m)+\widetilde J^{\ssup{q}}(p)-\langle q,\log g\rangle\Big\},\qquad p\in\Mcal_1(\R^d),
\end{aligned}
$$
where $\widetilde J^{\ssup{q}}$ is introduced in \eqref{ratefunctioncontract}, and we recall that $q=\mu_{0,\beta}$ if $I^{\ssup q}(\mu)<\infty$. As we mentioned below Theorem~\ref{main2}, we will show in Section~\ref{sec-proofYN} that $\widetilde J^{\ssup{q}}=J^{\ssup{q}}$. This finishes the proof.
\qed
\end{proofsect}

In the situation of Proposition~\ref{cormain}, the measure $\P_{\m,N}^{\ssup{{\rm sym}}}$ is not necessarily normalised, and no simple formula for its total mass seems available in general. Therefore, the following consequence of Proposition~\ref{cormain} seems helpful. When applied to discrete measures $\m$, it may have also some interesting consequences for related combinatorial questions. We also add a standard consequence of a large deviations principle: an identification of the minimiser of the rate function and a law of large numbers.

\begin{cor}\label{LLN} Under the assumptions of Proposition~\ref{cormain}, the following hold.
\begin{enumerate}

\item[(i)]\begin{equation}\label{combinatorics}
\begin{aligned}
\lim_{N\to\infty}&\frac{1}{N}\log\Big( \frac{1}{N!}\sum_{\s\in\Sym_N}\int_{(\R^d)^N}\prod_{i=1}^N \m(\d x_i)\,\prod_{i=1}^Ng(x_i,x_{\s(i)})\Big)\\
&=-\inf_{q\in\Mcal_1^{\ssup{\rm s}}(\R^d\times\R^d)}\Big\{H(q|\overline{q}\otimes\m)-\langle q,\log g\rangle\Big\}.
\end{aligned}
\end{equation}

\item[(ii)] The unique minimiser of the rate function $\mu\mapsto I^{\ssup{\rm sym}}_{\m}(\mu)-\langle \mu_{0,\beta},\log g\rangle$ is given by
\begin{equation}\label{mu*ident}
\mu^*=\int_{\R^d} \int_{\R^d}q^*(\d x,\d y)\,\P_{x,y}^{\beta}\circ B^{-1},
\end{equation} 
where $q^*\in\Mcal_1^{\ssup{\rm s}}$ is the unique minimiser of the formula on the right hand side of \eqref{combinatorics}.

\item[(iii)] {\em Law of large numbers:} Under the measure $ \P^{\ssup{\rm sym}}_{\m,N} $, normalised to a probability measure, the sequence $(L_N)_{N\in\N}$ converges in distribution to the measure $\mu^*$ defined in \eqref{mu*ident}.

\end{enumerate}
\end{cor}

\begin{proofsect}{Proof} According to Proposition~\ref{cormain}, the left hand side of \eqref{combinatorics} is equal to $-\inf_{\mu\in\Mcal_1(\Ccal)}(I^{\ssup{\rm sym}}_{\m}(\mu)-\langle \mu_{0,\beta},\log g\rangle)$ (use the large deviation principle for the measure of the event $\{L_N\in\Mcal_1(\Ccal)\}$). Use \eqref{Isymident} and substitute $q=\mu_{0,\beta}$ (recall that $I^{\ssup q}(\mu)=\infty$ otherwise) to see that this is equal to
\begin{equation}\label{infIsymg}
-\inf_{q\in\Mcal_1^{\ssup{\rm s}}(\R^d\times\R^d)}\Big[\Big\{H(q|\overline{q}\otimes\m)-\langle q,\log g\rangle\Big\}+\inf_{\mu\in\Mcal_1(\Ccal)\colon q=\mu_{0,\beta}}\sup_{\Phi\in\Ccal_{\rm b}(\Ccal)}\Big(\big\langle \mu,\Phi-\log \E_{\pi_0,\pi_\beta}^\beta\big[{\rm e}^{\Phi(B)}\big]\big\rangle\Big)\Big].
\end{equation}
It is easy to see that the latter infimum over $\mu$ is equal to zero. Indeed, for any $\mu$ pick $\Phi=0$ to see that \lq$\geq$\rq\ holds, and the choice $\mu=\int_{\R^d}\int_{\R^d}q(\d x,\d y)\,\P_{x,y}^\beta\circ B^{-1}$ and an application of Jensen's inequality shows that \lq$\leq$\rq\ holds. 

For proving (ii) and (iii) simultaneously, it suffices to show that $\mu^*$ is the unique minimiser of the rate function, $\mu\mapsto I^{\ssup{\rm sym}}_{\m}(\mu)-\langle \mu_{0,\beta},\log g\rangle $. Assume that $\mu$ is a zero of $ I^{\ssup{\rm sym}}_{\m} $. Since the map $q\mapsto H(q|\overline q\otimes \m)$ is known to have compact level sets, there is a $q^*\in\Mcal_1^{\ssup{\rm s}}(\R^d\times\R^d) $ that minimises the formula on the right hand side of \eqref{combinatorics}. 
Since in particular $I^{\ssup {q^*}}(\mu)<\infty$, we have $\mu_{0,\beta}=q^*$ and therefore
$$
0=I^{\ssup q^*}(\mu)=\sup_{\Phi\in\Ccal_{\rm b}(\Ccal)}\Big\langle\mu,\Phi-\log\E_{\pi_0,\pi_\beta}^{\beta}\big[{\rm e}^{\Phi(B)}\big]\Big\rangle.
$$
Hence, $\Phi\equiv 0$ is optimal in this formula. The variational equations yield, for any $h\in\Ccal_{\rm b}(\Ccal)$,
$$
\langle\mu,h\rangle=\big\langle\mu,\E_{\pi_0,\pi_\beta}^{\beta}[h(B)]\big\rangle.
$$
This identifies $\mu$ as $\mu^*$.
\qed
\end{proofsect}

There is an interesting by-product of Corollary~\ref{LLN} for the special case $g\equiv 1$, in which case it is easy to see that $q^*=\m\otimes\m$: In spite of strong correlations for fixed $N$ under $ \P^{\ssup{\rm sym}}_{\m,N} $, the initial and terminal locations $ B^{\ssup{1}}_0$ and $B^{\ssup{1}}_{\beta} $ of the first motion become independent in the limit $N\to\infty$.  One can prove this also in an elementary way, and also the fact that, for any $ k\in \N $ and for all $ i_1<i_2<\cdots<i_k $, the Brownian motions $ B^{\ssup{i_1}},\ldots,B^{\ssup{i_k}} $ under $ \P^{\ssup{\rm sym}}_{\m,N} $ become independent in the limit $ N\to\infty $.

\subsection{An important special case.}\label{sec-spcase}

%
%

\noindent In this section we consider an important special case of Proposition~\ref{cormain}(ii) that will be important for the applications in physics in Section~\ref{sec-Quantum}. We pick a large bounded closed box $\Lambda\subset\R^d$ and put $\m$ equal to $\Leb_\Lambda$,  the Lebesgue measure on $\Lambda$. Furthermore, we choose the function $g$ in Proposition~\ref{cormain} equal to $1/p_\beta(x,y)$, where $p_\beta$ is the Gaussian density in \eqref{Gaussian}. In other words, we replace the normalised Brownian bridge measure by the canonical, non-normalised one, $\mu^\beta_{x,y}$, introduced in \eqref{nnBBM}. That is, we look at the symmetrised measure
\begin{equation}\label{musymdef}
\mu_{\Lambda,N}^{\ssup{{\rm sym}}}=\frac{1}{N!}\sum_{\s\in\Sym_N}\;\int_{\Lambda}\cdots \int_{\Lambda} \d x_1\cdots\d x_N\,\bigotimes_{i=1}^N\mu^{\beta}_{x_i,x_{\s(i)}}.
\end{equation}
Apart from questions motivated from physics (see Section~\ref{sec-Quantum}), this measure is also mathematically interesting, see the discussion at the end of the present section. According to Proposition~\ref{cormain}(ii), the distribution of the mean of the normalised occupation measures, $Y_N$, under $\mu_{\Lambda,N}^{\ssup{{\rm sym}}}$ satisfies a large deviations principle (even though the term \lq distribution\rq\ is wrong since $\mu_{\Lambda,N}^{\ssup{{\rm sym}}}$ is not normalised). That is, we have
\begin{equation}\label{LDPmusym}
\lim_{N\to\infty}\frac 1N\log\big( \mu_{\Lambda,N}^{\ssup{{\rm sym}}}\circ Y_N^{-1}(\cdot)\big)=-\inf_{p\in \,\cdot} J^{\ssup{\rm sym}}_{\L}(p),
\end{equation}
in the weak sense on subsets of $\Mcal_1(\R^d)$, where we introduced
\begin{equation}\label{JsymLambdadef}
J^{\ssup{\rm sym}}_{\L}(p)=\inf_{q\in\Mcal_1^{\ssup{\rm s}}(\R^d\times\R^d)}\Big\{H(q|\overline{q}\otimes\Leb_\L)+ J^{\ssup{q}}_{p_\beta}(p)\Big\},\qquad p\in\Mcal_1(\R^d),
\end{equation}
and
\begin{equation}\label{JqLambdadef}
J^{\ssup{q}}_{p_\beta}(p)=\sup_{f\in\Ccal_{\rm b}(\R^d)}\Bigl\{\beta\langle f,p\rangle-\int_{\R^d}\int_{\R^d} q(\d x,\d y)\,\log\E_x\big[ {\rm e}^{\int_0^{\beta}f(B_s)\,\d s};B_\beta\in\d y\big]\big/\d y\Bigr\}.
\end{equation}
(In the notation of Section~\ref{sec-remarks}, $J^{\ssup{\rm sym}}_{\L}=J^{\ssup{\rm sym}}_{\Leb_\L,1/p_\beta}$.)
The main goal of the present section is to identify $J^{\ssup{\rm sym}}_{\L}$ in much easier and more familiar terms.
It turns out that $J^{\ssup{\rm sym}}_{\L}(p)$ is identical to the energy of the square root of the density of $p$, in the jargon of large deviations theory also sometimes called the {\it Donsker-Varadhan rate function}, $I_\L\colon \Mcal_1(\R^d)\to[0,\infty]$ defined by
\begin{equation}\label{Idef}
I_\L(p)=\begin{cases}\big\|\nabla \sqrt{\frac{\d p}{\d x}}\big\|_2^2,&\mbox{if }p\mbox{ has a density with square root in }H^1_0(\L^{\circ}),\\
\infty&\mbox{otherwise.}
\end{cases}
\end{equation}

\begin{theorem}\label{Jident} Let $ \L\subset\R^d $ be a bounded closed box. Then $ J^{\ssup{\rm sym}}_{\L}(p)= \beta I_{\L}(p)$ for any $p\in\Mcal_1(\R^d)$.
\end{theorem}

The proof of Theorem~\ref{Jident} is in Section~\ref{sec-JLambdaident}. In the theory and applications of large deviations, $I_\L$ plays an important role as the rate function for the normalised occupation measure of one Brownian motion (or, one Brownian bridge) in $\L$, in the limit as time to infinity \cite{DV75-83}, \cite{Ga77}. It is remarkable that, in Theorem~\ref{Jident}, in conjunction with Proposition~\ref{cormain}(ii), this function turns out also to govern the large deviations for the mean of the normalised occupation measures under the symmetrised measure $\mu_{\Lambda,N}^{\ssup{{\rm sym}}}$, in the limit of large number of motions. Let us give an informal discussion and interpretation of this fact. 

The measure $\mu_{\Lambda,N}^{\ssup{{\rm sym}}}$ in \eqref{musymdef} admits a representation which goes back to Feynman 1953 \cite{F53} and which we want to briefly discuss. 
Every permutation $\sigma\in\Sym_N$ can be written as a concatenation of cycles. Given a cycle $(i,\s(i),\s^2(i),\dots,\s^{k-1}(i))$ with $\s^k(i)=i$ and precisely $k$ distinct indices, the contribution coming from this cycle is independent of all the other indices. Furthermore, by the fact that $\mu^{\beta}_{x_i,x_{\s(i)}}$ is the conditional distribution given that the motion ends in $x_{\s(i)}$, this contribution (also executing the $k$ integrals over $x_{\s^l(i)}\in\Lambda$ for $l=k-1,k-2,\dots,0$) turns the corresponding $k$ Brownian bridges of length $\beta$ into one Brownian bridge of length $k\beta$, starting and ending in the same point $x_i\in\Lambda$ and visiting $\L$ at the times $\beta,2\beta,\dots,(k-1)\beta$. Hence, 
$$
\mu_{\Lambda,N}^{\ssup{{\rm sym}}}=\frac{1}{N!}\sum_{\s\in\Sym_N}\;\bigotimes_{k\in \N}\Big(\int_\Lambda \d y_k\,\mu^{k,\beta,\L}_{y_k,y_k}\Big)^{\otimes f_k(\sigma)},
$$
where $f_k(\s)$ denotes the number of cycles in $\s$ of length precisely equal to $k$, and $\mu^{k,\beta,\L}_{x,y}$ is the Brownian bridge measure $\mu^{k\beta}_{x,y}$ as in \eqref{nnBBM}, restricted to the event $\bigcap_{l=1}^k\{B_{l\beta}\in\L\}$. (See \cite[Lemma 2.1]{Gin71} for related combinatorial considerations.) If $f_N(\s)=1$ (i.e., if $\s$ is a cycle), then we are considering just one Brownian bridge $B$ of length $N\beta$, with uniform initial measure on $\L$, on the event $\bigcap_{l=1}^N\{B_{l\beta}\in\L\}$. Furthermore, $Y_N$ is equal to the normalised occupation measure of this motion. For such a $\s$, the limit $N\to\infty$ turns into a limit for diverging time, and the corresponding large-deviation principle of Donsker and Varadhan formally applies. 

If a permutation $\s$ does not contain a cycle of length $\approx N$, presumably, its contribution is quantified with a different rate. In this way, Theorem~\ref{Jident} says that the large-$N$ behaviour of $\mu_{\Lambda,N}^{\ssup{{\rm sym}}}\circ Y_N^{-1}$ is predominantly determined by all those permutations who consist of just one cycle of length $N$.

\subsection{Relation to quantum physics.}\label{sec-Quantum}

\noindent Let us now describe the relation of our work with the canonical ensemble of large systems of Bosons at positive temperature. We consider a system of $N$ non-interacting Bosons in a trap potential $ W $. The system is described by the Hamilton operator
\begin{equation}\label{Hdef}
\Hcal_N=\sum_{i=1}^N \big(-\Delta_i+W(x_i)\big),\qquad x_1,\dots,x_N\in\R^d,
\end{equation}
where the $i$-th Laplace operator, $ \Delta_i$, represents the kinetic energy of the $i$-th particle, and $W\colon\R^d\to[0,\infty]$ is the trap potential. The trace of the operator ${\rm e}^{-\beta \Hcal_N}$ is the canonical partition sum of the system at temperature $1/\beta$. However, the characteristic property of Bosons is expressed by the symmetry of any $N$-particle wave function under permutation of the coordinates. This in turn means that the partition sum of a system of $N$ Bosons is given by the trace of the restriction of ${\rm e}^{-\beta \Hcal_N}$ to the subspace of symmetric wave functions, denoted by ${\rm Tr}_+({\rm e}^{-\beta \Hcal_N})$. Via the Feynman-Kac formula, this trace is given as
\begin{equation}\label{trace}
\begin{aligned}
{\rm Tr}_+\big({\rm e}^{-\beta \Hcal_N}\big)&=\frac 1{N!}\sum_{\sigma\in\Sym_N}\int_{\R^{dN}}\d x_1\cdots \d x_N\int_{\Ccal^N}\Big(\bigotimes_{i=1}^N \mu^{\beta}_{x_i,x_{\sigma(i)}}\Big)(\d \omega)\,\exp\Big\{-\sum_{i=1}^N\int_0^\beta W(\omega_s^{\ssup i})\,\d s\Big\},
\end{aligned}
\end{equation}
where the canonical Brownian bridge measure was introduced in \eqref{nnBBM}, and we wrote $\omega=(\omega^{\ssup 1},\dots,\omega^{\ssup N})$.
One of our main results is an explicit formula for the logarithmic large-$N$ asymptotic of this trace, for a certain class of hard-wall traps $W$. This will be a consequence of Theorem~\ref{main2}, in conjunction with Theorem~\ref{Jident} and Varadhan's lemma. The main novelty of this result is the combined application of methods from combinatorics, variational analysis and the theory of large deviations to the study of the canonical ensemble.

Let us make some historical remarks. Feynman \cite{F53} analysed the partition function of an interacting Bose gas in terms of the statistical distribution of permutation cycles of particles and emphasised the roles of long cycles at the transition point. These arguments were pursued further by Penrose and Onsager \cite{PO56}. The arguments for the role of the cycle statistics have been known for a long time in various contexts, e.g., Ginibre \cite{Gin71} used them for virial expansion for quantum gases, Cornu \cite{C96} for the Mayer expansions for quantum Coulomb gases, and Ceperley \cite{C95} for numerical simulations for Helium via path integrals. In a couple of papers in the 1960ies, Ginibre studied the {\it grandcanonical\/} ensemble, where $N$ is a Poisson random variable; see the summary in \cite{Gin71}. His main interest was in \lq hard-wall\rq\ traps $W=\infty\1_{\L^{\rm c}}$, where $\L$ is a box, in the limit $\L\uparrow \R^d$. This corresponds in the canonical ensemble to the {\it thermodynamic limit}, i.e.,  the limit $N\to\infty$ with a box $\L=\L_N\uparrow \R^d$ coupled with $N$ in a way that the particle density per volume, $N/|\L_N|$, converges in $(0,\infty)$. A precise mathematical and quantitative formulation of the relation between Bose condensate and long cycles appeared only recently in work of S\"ut\"o \cite{S93}, \cite{S02} dealing with the ideal and mean field Bose gas. However, his methods are only applicable to the ideal gas or the mean field model, due to the difficult combinatorics of the cycle statistics.

The present paper introduces an alternative approach to a deeper understanding of the effect of the symmetrisation in large Boson systems by combining techniques coming from combinatorics and the theories of stochastic processes and large deviations. In future work, we will extend the techniques of the present paper to handle also boxes $\L=\L_N$ increasing to $\R^d$ as $N\to\infty$, as well as {interacting} Boson systems.

Let us return to the the symmetrised trace of $\Hcal_N$ in \eqref{trace}. We identify its large-$N$ asymptotic as follows. Given a box $\L\subset\R^d$, we denote by
\begin{equation}\label{eigenvdefi}
\l_\L(f)=\sup_{\varphi\in \Ccal^\infty(\R^d)\colon \supp(\varphi)\subset\L, \|\varphi\|_2=1}\Big(\langle f,\varphi^2\rangle-\|\nabla\varphi\|_2^2\Big)
\end{equation}
the principal (i.e., largest) eigenvalue of $\Delta+f$ in $\L$ with Dirichlet boundary condition.

\begin{theorem}\label{LDPW}
Fix $\beta\in(0,\infty)$ and let $\Lambda\subset\R^d$ be a bounded closed box. Let $W\colon\R^d\to\R\cup\{\infty\}$ be continuous in $\Lambda$ and equal to $\infty$ in $\Lambda^{\rm c}$. Then
\begin{equation}
\lim_{N\to\infty}\frac 1N \log\Big({\rm Tr}_+\big({\rm e}^{-\beta \Hcal_N}\big)\Big)=\beta\lambda_\L(W).
\end{equation}
\end{theorem}

\begin{proofsect}{Proof} Recall the measure $\mu_{\Lambda,N}^{\ssup{{\rm sym}}}$ from \eqref{musymdef} and the mean of the occupation measures, $Y_N$, from \eqref{YNdef}.
From \eqref{trace} we have that 
$$
{\rm Tr}_+\big({\rm e}^{-\beta \Hcal_N}\big)=\int_{\Ccal^N} {\rm e}^{-N\beta\langle W, Y_N\rangle}\1\{\supp(Y_N)\subset \L\}\,\d\mu^{\ssup{\rm sym}}_{\L,N}.
$$
Recall the large deviation principle of \eqref{LDPmusym}. Since the map $\Mcal_1(\L)\ni Y\mapsto \langle W,Y\rangle$ is bounded and continuous, we may apply Varadhan's lemma to deduce that
$$
\lim_{N\to\infty}\frac 1N \log\Big({\rm Tr}_+\big({\rm e}^{-\beta \Hcal_N}\big)\Big)=-\inf_{p\in \Mcal_1(\L)}\Big(J^{\ssup{\rm sym}}_{\L}(p)+\beta\langle W,p\rangle\Big).
$$
By Theorem~\ref{Jident}, we may replace $J^{\ssup{\rm sym}}_{\L}(p)$ by $\beta I_\L(p)$ defined in \eqref{Idef}. This gives that the right hand side is equal to $\beta\sup_{p\in \Mcal_1(\L)}[\langle W,p\rangle-I_\L(p)]$. The substitution $\varphi^2(x)\,\d x=p(\d x)$ and a glance at \eqref{eigenvdefi} yield that this is equal to $\beta\l_\L(W)$.
\qed
\end{proofsect}

\section{Identification of $J^{\ssup{\rm sym}}_{\L}$}\label{sec-JLambdaident}

In this section, we prove Theorem~\ref{Jident}.
First we consider the rate function $J^{\ssup{\rm sym}}_{\L}$ defined in \eqref{JsymLambdadef}. By $B_{[0,\beta]}=\{B_s\colon s\in[0,\beta]\}$ we denote the path of the Brownian motion $B$, and $\Ccal(\Lambda)$ denotes the set of continuous functions $\L\to\R$.

\begin{lemma}\label{Jrestr} Fix $\beta\in(0,\infty)$ and a bounded closed box $\L\subset\R^d$. Then, for all $ p\in\Mcal_1(\R^d) $ having support in $\L$,
\begin{equation}\label{JqLambdaident}
J^{\ssup{\rm sym}}_{\L}(p)=
\inf_{q\in\Mcal_1^{\ssup{\rm s}}(\L\times\L)}\Big\{H(q|\overline{q}\otimes\Leb_\L)+J^{\ssup q}_{\L,p_\beta}(p)\Big\},
\end{equation}
where
$$
J^{\ssup q}_{\L,p_\beta}(p)=\sup_{f\in\Ccal(\L)}\Bigl\{\beta\langle f,p\rangle-\int_{\L}\int_{\L} q(\d x,\d y)\,\log\E_x\big[ {\rm e}^{\int_0^{\beta}f(B_s)\,\d s}\1\{B_{[0,\beta]}\subset\L\};B_\beta\in\d y\big]\big/\d y\Bigr\}
$$
\end{lemma}

\begin{proofsect}{Proof} Note that $H(q|\overline{q}\otimes\Leb_\L)=\infty$ if the support of $q$ is not contained in $\L\times\L$. Hence, in \eqref{JsymLambdadef} we need to take the infimum over $q$ only on the set $\Mcal_1^{\ssup{\rm s}}(\L\times\L)$.
From an inspection of the right hand side of \eqref{JqLambdadef} it follows that the function $f$ in the supremum must be taken arbitrarily negative outside $\L$ to approximate the supremum. Hence, we may add in the expectation the indicator on the event that the Brownian motion does not leave $\L$ by time $\beta$. But then the values of $f$ outside $\L$ do not contribute. This shows that we need to consider only continuous functions $f$ that are defined on $\L$; in other words,  \eqref{JqLambdaident} holds.
\qed
\end{proofsect}

\begin{proofsect}{Proof of Theorem~\ref{Jident}} We start from \eqref{JqLambdaident} and proceed in three steps: (1) we show that $ J^{\ssup{\rm sym}}_{\L}(p)\ge \beta I_{\L}(p)$ for any $ p\in\Mcal_1(\R^d) $ with support in $ \L$, (2) we show that the complementary inequality, $ J^{\ssup{\rm sym}}_{\L}(p)\le \beta I_\L(p)$, holds if $ \varphi:=\sqrt{\frac{\d p}{\d x}}$ exists in $ \Ccal^2(\R^d) $ with $ \frac{\Delta\varphi}{\varphi}\in\Ccal(\L) $, and (3) we approximate an arbitrary $p\in\Mcal_1(\R^d)$ satisfying $\varphi\in H_0^1(\L^\circ)$ with suitable smooth functions.

Let us turn to the details. For $ f\in\Ccal(\L) $, let $ \varphi_f$ be the unique positive $ L^2$-normalised eigenfunction of $ \Delta+f $ in $ L^2(\L) $ with zero boundary condition and corresponding eigenvalue $\l_{\L}(f)$; see \eqref{eigenvdefi}. Then
\begin{equation}\label{densitygirsanov}
D_{\beta}^{\ssup{f}}:={\rm e}^{\int_0^{\beta}f(B_s)\,\d s}{\rm e}^{-\beta\l_{\L}(f)}\1\{B_{[0,\beta]}\subset\L\}\frac{\varphi_f(B_{\beta})}{\varphi_f(B_{0})}
\end{equation} 
defines a martingale $ (D_{\beta}^{\ssup{f}})_{\beta\ge 0} $ under $\P_x$ for any $x\in\L$ with respect to the canonical Brownian filtration (see \cite[Prop. VIII.3.1]{RY99}). Substituting $D_{\beta}^{\ssup{f}}$ on the right hand side of \eqref{JqLambdaident} and using the marginal property of $ q $ (i.e., $\int\int q(\d x,\d y)\log(\varphi(y)/\varphi(x))=0$), we see that
$$
J^{\ssup{q}}_{\L,p_\beta}(p)=\sup_{f\in\Ccal(\L)}\Bigl(\beta[\langle f,p\rangle - \l_{\L}(f)]-\int_{\L}\int_{\L}q(\d x,\d y)\,\log\frac{\E_{x}\big[D_{\beta}^{\ssup{f}};B_{\beta}\in \d y\big]}{\d y}\Bigr),
$$ 
where $\E_x$ denotes expectation with respect to a Brownian motion with generator $\Delta$ starting at $x$. Substituting this in \eqref{JqLambdaident}, we obtain that
\begin{equation}\label{infqsupf}
J^{\ssup{\rm sym}}_{\L}(p)=\inf_{q\in\Mcal_1^{\ssup{\rm s}}(\L\times\L)}\sup_{f\in\Ccal(\L)}\Bigl(\beta[\langle f,p\rangle - \l_{\L}(f)]+\int_{\L}\int_{\L}q(\d x,\d y)\,\log\frac{q(\d x,\d y)}{\overline{q}(\d x)\E_{x}\big[D_{\beta}^{\ssup{f}};B_{\beta}\in \d y\big]}\Bigr).
\end{equation}
By the martingale property of $ (D_{\beta}^{\ssup{f}})_{\beta\ge 0} $, the measure $ \E_x[D_{\beta}^{\ssup{f}};B_{\beta}\in\d y] $ is a probability measure on $ \L $ for any $ x \in\L $. Hence, the double integral in \eqref{infqsupf} is an entropy between probability measures and therefore nonnegative, by Jensen's inequality. This shows that
\begin{equation}\label{JsymLambdaident}
J^{\ssup{\rm sym}}_{\L}(p)\ge \beta \sup_{f\in\Ccal(\L)}\Bigl(\langle f,p\rangle-\l_{\L}(f)\Bigr).
\end{equation}
Note that the map $f\mapsto \l_{\L}(f) $ is the Legendre-Fenchel transform of  $ I_\L $, as is seen from the Rayleigh-Ritz principle in \eqref{eigenvdefi}. According to the Duality Lemma \cite[Lemma~4.5.8]{DZ98}, the r.h.s.\ of \eqref{JsymLambdaident} is therefore equal to $ \beta I_\L(p) $ since it is equal to the Legendre-Fenchel transform of $\l_{\L}$. Hence, we have shown that $ J^{\ssup{\rm sym}}_{\L}(p)\ge \beta I(p)$ for any $ p\in\Mcal_1(\R^d) $ with support in $ \L$.

Now we proceed with the second step. Let $ \varphi=\sqrt{\frac{\d p}{\d x}} $ be in $ \Ccal^2(\R^d) $ such that $ f^*=-\frac{\Delta\varphi}{\varphi} $ is in $ \Ccal(\L) $. Then $ (\Delta+f^*)\varphi=0 $ in $ \L $. In other words, $ \varphi=\varphi_{f^*} $ is the unique positive normalised eigenfunction of $\Delta+f^*$ in $\L$ with corresponding eigenvalue $ \l_{\L}(f^*)=0 $. Consider the measure
$$
\begin{aligned}
q^*(\d x,\d y)&=\varphi(x)\varphi(y)\int_\Ccal{\rm e}^{\int_0^{\beta}f^*(\omega_s)\,\d s}\1\{\omega_{[0,\beta]}\subset\L\}\,\mu_{x,y}^\beta(\d \omega)\, \d x\d y\\
&=\varphi(x)\varphi(y)\E_x\Big[{\rm e}^{\int_0^{\beta}f^*(B_s)\,\d s}\1\{B_{[0,\beta]}\subset\L\}\,; B_\beta\in \d y\Big]\,\d x
\end{aligned}
$$ 
on $ \L\times\L $. Then $q^*$ is obviously symmetric. With the help of the martingal property of $(D_{\beta}^{\ssup{f}})_{\beta\ge 0} $ , the marginal of $q^*$ is identified as
$$
\overline{q^*}(\d x)=\varphi(x)\E_{x}\Big[{\rm e}^{\int_0^{\beta}f^*(B_s)\,\d s}\1\{B_{[0,\beta]}\subset\L\}\varphi(B_{\beta})\Big]\,\d x=\varphi^2(x)\,\d x=p(\d x).
$$ 
Hence $ q^*\in\Mcal_1^{\ssup{\rm s}}(\L\times\L) $. Using this $ q^* $ in \eqref{infqsupf}, we obtain, also using  the marginal property of $ q^* $,
\begin{equation}\label{varupper}
J^{\ssup{\rm sym}}_{\L}(p)\le \sup_{f\in\Ccal(\L)}\Bigl(\beta\langle f,p\rangle+\int_{\L}\int_{\L}q^*(\d x,\d y)\,\log\frac{ \E_{x}\big[{\rm e}^{\int_0^{\beta}f^*(B_s)\,\d s}\1\{B_{[0,\beta]}\subset\L\}\,;B_\beta\in\d y\big]}{\E_{x}\big[{\rm e}^{\int_0^{\beta}f(B_s)\,\d s}\1\{B_{[0,\beta]}\subset\L\}\,;B_\beta\in\d y\big]}\Bigr).
\end{equation}
Now we show that the variational problem on the r.h.s.~of \eqref{varupper} is solved precisely in $ f=f^* $. Indeed,  by strict concavity we only have to show that $f=f^*$ solves the variational equation, which reads
$$
\forall\, h\in\Ccal(\L)\colon\qquad \beta\langle h,p\rangle=\int_{\L}\int_{\L}q^*(\d x,\d y)\,\frac{ \E_{x}\Big[\Big(\int_0^{\beta}h(B_s)\,\d s\Big){\rm e}^{\int_0^{\beta}f(B_s)\,\d s}\1\{B_{[0,\beta]}\subset\L\}\,;B_\beta\in\d y\Big]}
{\E_{x}\Big[{\rm e}^{\int_0^{\beta}f(B_s)\,\d s}\1\{B_{[0,\beta]}\subset\L\}\,;B_\beta\in\d y\Big]}.
$$ 
This is indeed solved for $ f=f^* $, because the right hand side of the variational equation then equals
$$
\begin{aligned}
\int_{\L}\int_{\L}\d x \d y\,&\varphi(x)\varphi(y)\E_{x}\Big[\Big(\int_0^{\beta}h(B_s)\,\d s\Big){\rm e}^{\int_0^{\beta}f^*(B_s)\,\d s}\1\{B_{[0,\beta]}\subset\L\}\,;B_\beta\in\d y\Big]\\
&=\int_{\L}\d x\,\varphi^2(x)\E_x\Big[D_{\beta}^{\ssup{f^*}}\int_0^{\beta}h(B_s)\,\d s\Big]
=\int_0^{\beta}\d s\,\widehat \E^{\ssup{f^*}}\big[h(B_s)\big],
\end{aligned}
$$ 
where $ \widehat\E^{\ssup{f^*}} $ is expectation with respect to the Girsanov transform with martingale $ (D_{\beta}^{\ssup{f^*}})_{\beta\ge 0} $ defined in \eqref{densitygirsanov}, starting in its invariant distribution, $ \varphi^2(x)\,\d x=p(\d x) $. Note that the transformed Brownian motion does not leave $\L$ and is stationary, when started with distribution $p$. Hence, $ \widehat\E^{\ssup{f^*}}[h(B_s)]=\langle h,p\rangle $ for any $ s\in [0,\infty) $. Therefore $f=f^* $ solves the variational equation and is a maximiser on the right hand side of \eqref{varupper}. Hence, 
$$
J^{\ssup{\rm sym}}_{\L}(p)\le \beta\langle f^*,p\rangle=-\beta\langle\Delta\varphi,\varphi\rangle=\beta||\nabla\varphi||_2^2=\beta I_\L(p).
$$

Now we finish the proof. Let $p\in\Mcal_1(\R^d)$ be arbitrary with support in $\L$. We need to show that $J^{\ssup{\rm sym}}_{\L}(p)\le\beta I_\L(p)$. Certainly, we may assume that $\varphi=\sqrt {\frac{\d p}{\d x}}$ exists and lies in $H_0^1(\L^\circ)$. Hence, there is a sequence of smooth functions $\varphi_n\in\Ccal^\infty$ such that $\supp(\varphi_n)\subset \L^\circ$ and $\varphi_n\to \varphi$ in $H^1$-norm. 

We need to approximate $\varphi_n$ with suitable smooth functions $\widetilde \varphi_n$ such that $\frac{\Delta\widetilde \varphi_n}{\widetilde \varphi_n}$ is continuous in $\L$. For this purpose, choose $\delta_n>0$ such that $\supp(\varphi_n)\subset \L_{\delta_n}=\{x\in\L\colon \dist(x,\L^{\rm c})\geq \delta_n\}$. Pick some small $\eps_n>0$ and some smooth function $\kappa_n\colon\R^d\to[0,1]$ satisfying $\supp(\kappa_n)\subset \L^\circ$ and $\kappa_n(x)=1$ for $x\in \L_{\delta_n}$ such that $\Delta\kappa_n/\kappa_n$ and $\nabla\kappa_n/\kappa_n$ are continuous in $\L$. Then we put $\widetilde \varphi_n=\kappa_n^{1/2}(\varphi_n^2+\eps_n)^{1/2}$. Then $\widetilde \varphi_n$ is smooth with support in $\L^\circ$, and $\widetilde \varphi_n$ converges towards $\varphi$ in $L^2$. Furthermore, $\frac{\Delta\widetilde \varphi_n}{\widetilde \varphi_n}$ is continuous in $\L$. If $\eps_n$ is small enough (depending on $\delta_n$ and $\kappa_n $ only),  we also have that $\limsup_{n\to\infty} \|\nabla\widetilde \varphi_n\|_2\leq \|\nabla\varphi\|_2$ (use that $\varphi_n(x)=0$ for $x\in\L\setminus  \L_{\delta_n}$).

Along a suitable subsequence, $\widetilde\varphi_n$ converges almost everywhere to $\varphi$. Let $p_n\in\Mcal_1(\R^d)$ be the measures with density $\widetilde\varphi_n/\|\widetilde\varphi_n\|_2$. By the second step of the proof, we have $J^{\ssup{\rm sym}}_{\L}(p_n)\leq\beta I_\L(p_n)$ for any $n\in\N$. With the help of Fatou's lemma, we see that $\langle f, p\rangle\leq \liminf_{n\to\infty}\langle f, p_n\rangle$ for any $f\in\Ccal(\L)$. Hence, it is clear that
$$
J^{\ssup{\rm sym}}_{\L}(p)\leq \liminf_{n\to\infty}J^{\ssup{\rm sym}}_{\L}(p_n)\leq \liminf_{n\to\infty}\beta I_\L(p_n)\leq\beta I_\L(p).
$$
This ends the proof.
\end{proofsect}
\qed

\section{Proof of Theorem~\ref{main}}\label{sec-mainproof}

In this section we prove Theorem \ref{main}. The following is a reformulation of that theorem, making explicit what a large deviations principle is.

\begin{theorem}[Reformulation of Theorem \ref{main}]\label{mainnew} Fix $\beta\in(0,\infty)$ and $\m\in\Mcal_1(\R^d)$.

\begin{enumerate}

\item[(i)] Assume that $ \m $ has compact support. Then, for any open set $G\subset\Mcal_1(\Ccal) $,
\begin{equation}\label{lowbound}
\liminf_{N\to\infty}\frac 1N\log\P^{\ssup{\rm sym}}_{\m,N}(L_N\in G)\geq -\inf_{\mu\in G}I^{\ssup{\rm sym}}_{\m}(\mu).
\end{equation}

\item[(ii)] For any compact set $F\subset\Mcal_1(\Ccal) $,
\begin{equation}\label{uppbound}
\limsup_{N\to\infty}\frac 1N\log\P^{\ssup{\rm sym}}_{\m,N}(L_N\in F)\leq -\inf_{\mu\in F}I^{\ssup{\rm sym}}_{\m}(\mu),
\end{equation}

\item[(iii)] The sequence of distributions of $L_N$ under $\P^{\ssup{\rm sym}}_{\m,N}$ is exponentially tight.
\end{enumerate}
\end{theorem}

An outline of the proof is in Section~\ref{sec-outline}. The respective parts of Theorem~\ref{mainnew} are proved in the remaining subsections.

\subsection{Outline of the proof.}\label{sec-outline}

\noindent Let us briefly outline the main idea in the proof of Theorem~\ref{mainnew}(i) and (ii). The methods of the proof consist of a discretisation argument similar to \cite{KM02}, combined with combinatorial considerations (see, e.g.~\cite{A01}) and large-deviations arguments.

For technical reasons, we first replace $\R^d$ by a large ball $\Lambda$, which contains $\supp(\m)$ in the proof of the lower bound, respectively is later sent to $\R^d$ in the proof of the upper bound. The first main idea is that there is no problem in proving a large deviations principle for $L_N$ under a measure of the form $\bigotimes_{r,s}(\P^\beta_{x_r,x_s})^{\otimes N \eta(r,s)}$ if the integers $N \eta(r,s)$ sum up to one over a finite index set of $r$'s and $s$'s, and if the $x_r\in \Lambda$ are fixed. Such a large deviations principle follows in a standard way from the G\"artner-Ellis theorem \cite[Th.~4.5.20]{DZ98}. Here we consider the mean of $N$ random variables $\delta_{B^{\ssup i}}$, who are independent, but not identically distributed, but the number of distributions is fixed.

However, the problem is that, in \eqref{Psymdef}, for fixed $\s\in\Sym_N$ and for fixed integration variables $x_1,\dots,x_N$, the variety of measures $\P^\beta_{x_r,x_s}$ appearing is much too large for an application of this idea; the complexity is too high. Therefore, we introduce a partition of $\Lambda$ into finitely many small subsets $U_r$. For any $s$, we treat all $\P^\beta_{x_i,x_{\s(i)}}$ as equal if $x_{\s(i)}$ lies in the same $U_s$. More precisely, we relax the condition that the motion ends precisely in $x_{\s(i)}$ by the requirement that it ends somewhere in $U_s$. If the fineness is small enough, the replacement error will be small. Then we integrate out with respect to $\m$ over all $x_i$ within their partition set $U_r$, say. In this way, we have replaced the \lq microscopic\rq\ picture of the $\P^\beta_{x_i,x_{\s(i)}}$ by a \lq macroscopic\rq\ one that registers only the partition sets, $U_r$ and $U_s$, in which the motion starts and terminates, respectively. This we do for any $r$ and $s$ simultaneously. In this way, we now obtain a {\it finite\/} complexity of different types of Brownian bridges, ordered according to their initial and terminal partition sets. Say, for any $r,s$, we have $N\eta(r,s)$ motions that start in $U_r$ and end in $U_s$. Certainly, we have to sum over all admissible multi-indices $\eta$.

So far, we have not talked about the role of the permutations. For a given admissible $\eta$ and given integration variables $x_1,\dots,x_N$, only those $\s\in\Sym$ contribute that have the property that, for any $r,s$, precisely $N\eta(r,s)$ indices $i$ satisfy $x_i\in U_r$ and $x_{\s(i)}\in U_s$. The point is that these $\s$'s yield precisely the same contribution. Therefore only their cardinality is to be examined. This causes some combinatorial work that can be done in an elementary way. The result is expressible in terms of quotients of factorials, which can asymptotically be well approximated by entropies, using Stirling's formula. In this way, one arrives at a discrete version of the variational formula on the right hand sides of \eqref{lowbound} and \eqref{uppbound}, respectively. Some analytical work has to be done when letting the fineness of the partition vanish.

\subsection{Proof of Theorem~\ref{mainnew}(i).}\label{sec-lowproof}

\noindent We have to introduce some notation, which will be used frequently in the entire section. For any compact set $\Lambda\subset\R^d $ and any partition $\Ucal=\{U_r\colon r\in \Sigma\}$ of $\Lambda$ we denote by $f_\Ucal=\max_{r\in\Sigma}\diam(U_r)$ its fineness. We always tacitly assume that the sets $U_r$ are measurable and satisfy $\m(U_r)>0$ for any $r\in\Sigma$. In particular, for the proof of the lower bound \eqref{lowbound}, we choose $ \Lambda $ such that $ \supp\,(\m)\subset \Lambda $. By $m(r)=m_\Ucal(r)=\m(U_r)/\m(\Lambda)$ we denote the coarsened and normalised version of $\m$ on $\Lambda$; hence $m\in\Mcal_1(\Sigma)$. Furthermore, we introduce the set of probability measures $\eta$ on $ \Sigma^2 $ having equal marginals $ \overline{\eta}(r)=\sum_{s\in\Sigma}\eta(r,s) $:
\begin{equation}\label{pairmeasures}
\Mcal_1^{\ssup{\rm s}}(\Sigma^2)=\Bigl\{\eta\in\Mcal_1(\Sigma^2)\colon\sum_{s\in\Sigma}\eta(r,s)=\sum_{s\in\Sigma}\eta(s,r),\,\forall r\in\Sigma\Bigr\}.
\end{equation}
Additionally, let $ \Mcal_1^{\ssup{N}}(\Sigma^2)=\Mcal_1(\Sigma^2)\cap\frac{1}{N}\N_0^{\Sigma^2} $ be the set of those pair measures $\eta$ such that all the numbers $N\eta(r,s)$ are integers, and $ \Mcal_1^{\ssup{{\rm s},N}}(\Sigma^2)=\Mcal_1^{\ssup{\rm s}}(\Sigma^2)\cap\frac{1}{N}\N_0^{\Sigma^2} $.
We also need to introduce the probability measure 
\begin{equation}\label{low3}
P^{\beta}_{\eta,N,\Ucal}=\bigotimes_{r,s\in\Sigma}\Bigl(\P_{U_r,U_s}^{\beta}\Bigr)^{\otimes N\eta(r,s)},\qquad \eta\in \Mcal_1^{\ssup{{\rm s},N}}(\Sigma^2),
\end{equation} 
where
\begin{equation}\label{PUUdef}
\P_{U_r,U_s}^{\beta}=\int_{U_r}\frac{\m(\d x)}{\m(U_r)}\P_x(\,\cdot\,|B_{\beta}\in U_s)=\int_{U_r\times U_s}\frac{\m(\d x)p_{\beta}(x,y)\,\d y}{\m(U_r)\int_{U_s}p_{\beta}(x,z)\,\d z}\P_{x,y}^{\beta}\qquad\mbox{for } r,s\in \Sigma,
\end{equation}
is a coarsened version of the Brownian bridge measure, see \eqref{Gaussian}.
The entropy of measures on $\Sigma^2$ is also denoted by $H$, that is,
\begin{equation}\label{entropy}
H(\eta|\overline{\eta}\otimes m)=\sum_{r,s\in \Sigma}\eta(r,s)\log\frac{\eta(r,s)}{\overline{\eta}(r) m(s)}
\end{equation} 
is the relative entropy of $ \eta\in\Mcal_1^{\ssup{\rm s}}(\Sigma^2) $ with respect to the product $ \overline{\eta}\otimes m $. By $ \d $  we denote the L\'{e}vy metric on $\Mcal_1(\Ccal)$, which generates the weak topology; see \eqref{Levymetric} below. By $\dist(\mu,A)=\inf_{\nu\in A}\d(\mu,\nu)$ we denote the distance to a set $A\subset \Mcal_1(\Ccal) $.

Our first main step is presented now; it basically summarises all combinatorial arguments needed in the proof of Theorem~\ref{mainnew}(i).

\begin{prop}\label{firstproplower} Let an open set $ G\subset\Mcal_1(\Ccal) $ be given, fix $\delta>0$ and put $G_\delta=\{\mu\in G\colon  \dist(\mu,G^{\rm c})>\delta\}$. Let $\Lambda\subset\R^d $ be a compact set  such that $\supp\,(\m)\subset \Lambda $ and pick any partition $\Ucal=\{U_r\colon r\in \Sigma\}$ of $\Lambda$ with fineness $\leq \delta$. Then, for any $N\in\N$,
\begin{eqnarray}\label{summarising}
\begin{aligned}
\P^{\ssup{\rm sym}}_{\m,N}&(L_N\in G)
&\ge (CN)^{-\frac 12(\sharp\Sigma)^2}\sum_{\eta\in\Mcal_1^{\ssup{{\rm s},N}}(\Sigma^2)}{\rm e}^{-N H(\eta|\overline{\eta}\otimes m)}P^{\beta}_{\eta,N,\Ucal}(L_N\in G_{\delta}),
\end{aligned}
\end{eqnarray} 
where $C\in(0,\infty)$ is an absolute constant, introduced in \eqref{StirlingC} below.
\end{prop}

\begin{proofsect}{Proof}
According to \eqref{Psymdef},
\begin{equation}
\P^{\ssup{\rm sym}}_{\m,N}(L_N\in G)= \frac{1}{N!}\sum_{\s\in\Sym_N}\int_{\Lambda^N}\m(\d x_1)\cdots\m(\d x_N)\,\Bigl(\bigotimes_{i=1}^N\P_{x_i,x_{\s(i)}}^{\beta}\Bigr)(L_N\in G).
\end{equation}
Let $\Ucal=\{U_r\colon r\in\Sigma\}$ a partition of $\Lambda$. Hence, $\Lambda=\bigcup_{r\in\Sigma} U_r$, and $U_r\cap U_s=\emptyset $ for $r\not=s$. 
We split the integration over $\Lambda$ into sums of integrations over the subsets as
\begin{equation}\label{split}
\int_{\Lambda^N}=\sum_{r_1,\ldots,r_N\in\Sigma}\;\int_{U_{r_1}}\cdots\int_{U_{r_N}}.
\end{equation} 
For the ease of notation we sometimes write $r(i)$ instead of $r_i$. Using \eqref{split}, we can estimate
\begin{equation}\label{low1}
\begin{aligned}
\P^{\ssup{\rm sym}}_{\m,N}(L_N\in G)& =\frac{1}{N!}\sum_{\s\in\Sym_N}\sum_{r_1,\ldots,r_N\in \Sigma}\int_{U_{r_1}}\m(\d x_1)\cdots\int_{U_{r_N}}\m(\d x_N)\bigotimes_{i=1}^N\Bigl(\P_{x_{r(i)},x_{r(\s(i))}}^{\beta}\Bigr)(L_N\in G)\\
&\geq \frac{1}{N!}\sum_{\s\in\Sym_N}\sum_{r_1,\ldots,r_N\in\Sigma}\int_{U_{r_1}}\m(\d x_1)\cdots\int_{U_{r_N}}\m(\d x_N) \\
&\qquad\qquad \inf_{y_i\in U_{r(i)}, 1\le i\le N} \Bigl(\bigotimes_{i=1}^N\P_{x_{r(i)},y_{\sigma(i)}}^{\beta}\Bigr)(L_N\in G).
\end{aligned}
\end{equation} 
Introduce 
\begin{equation}\label{PUydef}
\P^{\beta}_{U_r,y}=\int_{U_r}\frac{\m(\d x)}{\m(U_r)}\P_{x,y}^{\beta},\qquad y\in\R^d, r\in\Sigma.
\end{equation} 
Using this notation, \eqref{low1} reads
\begin{eqnarray}\label{estintegralstart}
\begin{aligned}
\P^{\ssup{\rm sym}}_{\m,N}(L_N\in G)\ge \frac{1}{N!}\sum_{\s\in\Sym_N}\sum_{r_1,\ldots,r_N\in\Sigma}\prod_{i=1}^Nm(r_i)\inf_{y_i\in U_{r(i)}, 1\le i\le N}\Bigl(\bigotimes_{i=1}^N\P_{U_{r(i)},y_{\s(i)}}^{\beta}\Bigr)(L_N\in G).
\end{aligned}
\end{eqnarray}
In the following we replace the measures $ \P_{U_{r(i)},y_i}^{\beta} $ on the right hand side of (\ref{estintegralstart}) with the measures $\P_{U_r,U_s}^{\beta}$ defined in \eqref{PUUdef}.
To that end, we need to make the set $G$ a bit smaller, more precisely, we have to replace it by the set $G_\delta$.
Recall the L\'{e}vy metric $ \d $ on the Polish space $ \Mcal_1(\Ccal) $ \cite{DS01}, defined for any two probability measures $ \mu,\nu\in\Mcal_1(\Ccal) $ as
\begin{equation}\label{Levymetric}
\d(\mu,\nu)=\inf\{\delta >0\colon\mu(\Gamma)\le \nu(\Gamma^{\delta})+\delta\,\mbox{ and }\, \nu(\Gamma)\le\mu(\Gamma^{\delta})+\delta\,\mbox{ for all }\, \Gamma=\overline{\Gamma}\subset\Ccal\}, 
\end{equation} 
where $ F^{\delta}=\{\mu\in\Mcal_1(\Ccal)\colon\dist(\mu,F)\leq\delta\} $ is the closed $\delta$-neighbourhood of $ F $. 

\begin{lemma}\label{estimationsetendpoints}
Let $ \delta> 0 $. Pick a partition $\Ucal$ with fineness $ f_{\Ucal}\le\delta $. Then, for any $ r_1,\ldots,r_N\in\Sigma $, any $ y_1\in U_{r(1)},\ldots,y_N\in U_{r(N)} $ and any $ \s\in\Sym_N $, we have,\\
(i) for any open set $ G\subset\Mcal_1(\Ccal) $,
\begin{equation}
\bigotimes_{i=1}^N\P^{\beta}_{U_{r(i)},y_{\s(i)}}(L_N\in G)\ge \bigotimes_{i=1}^N\P^{\beta}_{U_{r(i)},U_{r(\s(i))}}(L_N\in G_{\delta}),
\end{equation}
\noindent (ii) for any closed set $ F\subset\Mcal_1(\Ccal) $,
\begin{equation}
\bigotimes_{i=1}^N\P^{\beta}_{U_{r(i)},y_{\s(i)}}(L_N\in F)\le \bigotimes_{i=1}^N\P^{\beta}_{U_{r(i)},U_{r(\s(i))}}(L_N\in F^{\delta}).
\end{equation}
\end{lemma}

\begin{proofsect}{Proof}
For any $ i=1,\ldots, N, $ we construct a Brownian bridge $ B^{\ssup{i}} $ under the measure $ \P_{U_{r(i)},y_{\s(i)}}^\beta $ and a conditioned Brownian motion $ \widetilde B^{\ssup{i}} $ under the measure $ \P_{U_{r(i)},U_{r(\s(i))}}^{\beta} $ jointly on one probability space as follows. Let $ W_1,\ldots, W_N $ be independent Brownian motions on $ [0,\beta] $ starting with distributions $ \m_{U_{r(1)}},\ldots,\m_{U_{r(N)}}, $ respectively, where $ \m_{U_r}:=\frac{\m|_{U_r}}{\m(U_r)}$ for $r\in\Sigma.$  Put, for $t\in [0,\beta]$,
\begin{eqnarray*}
B_t^{\ssup{i}}&=&W_i(t)+\frac{t}{\beta}(y_{\s(i)}-W_i(\beta)),\\
\widetilde B_t^{\ssup{i}}&=&W_i(t)+\frac{t}{\beta}(Z_{\s(i)}-W_i(\beta)),
\end{eqnarray*}
where $ Z_{\s(i)} $ has distribution  $ \P_{U_{r(i)},U_{r(\s(i))}}^\beta\circ B_{\beta}^{-1} $ and is independent of $ W_1,\ldots, W_N.$ Since $ \diam\; U_{r(\s(i))} \le \delta $, we have $ ||B^{\ssup{i}}-\widetilde B^{\ssup{i}}||_{\infty}\le \delta $. This implies that 
$$ 
\d\Bigl(\frac{1}{N}\sum_{i=1}^N\delta_{B^{\ssup{i}}},\frac{1}{N}\sum_{i=1}^N\delta_{\widetilde B^{\ssup{i}}}\Bigr)\le \delta,
$$ and therefore  both assertions.\qed
\end{proofsect}

Using Lemma \ref{estimationsetendpoints}, we arrived, for any $\delta>0$, at the estimate
\begin{eqnarray}\label{low2}
\begin{aligned}
\P^{\ssup{\rm sym}}_{\m,N}(L_N\in G)\ge \frac{1}{N!}\sum_{\s\in\Sym_N}\sum_{r_1,\ldots,r_N\in\Sigma}\prod_{i=1}^N m(r_i)\Bigl(\bigotimes_{i=1}^N\P_{U_{r_i},U_{r(\s(i))}}^{\beta}\Bigr)(L_N\in G_{\delta}),
\end{aligned}
\end{eqnarray}
for any partition $\Ucal$ of $\L$ with fineness $f_\Ucal\leq \delta$, where $G_{\delta}=\{\mu\in G\colon \dist (\mu,G^{\rm c}) >\delta\} $.

An important observation is that the probability term on the right of \eqref{low2} does not really depend on the full information contained in $\sigma$ and $r_1,\dots, r_N$, but only on the frequency of all the pairs $\in\Sigma^2$ in the sequence $(r_1,r(\sigma(1)),\dots,(r_N,r(\s(N)))$. In order to take advantage of this observation, we rewrite the right hand side of \eqref{low2} in terms of probability measures $\eta$ on $\Sigma^2$. For $ \eta\in \Mcal_1^{\ssup{{\rm s},N}}(\Sigma^2) $ and $ {R}=(r_1,\ldots,r_N)\in\Sigma^N $, let
\begin{equation}
\Sym_N({R},\eta)=\left\{\s\in\Sym_N\colon\sharp\{i\colon r_i=r, r_{\s(i)}=s\}=N\eta(r,s),\,\forall r,s\in\Sigma\right\}
\end{equation} 
be the set of those permutations $\s$ such that $\eta$ is equal to the empirical measure of the sequence $(r_1,r(\sigma(1)),\dots,(r_N,r(\s(N)))$. Recall the path probability measures $P^{\beta}_{\eta,N,\Ucal}$ defined in \eqref{low3} and note that 
\begin{equation}\label{low4}
\bigotimes_{i=1}^N\P_{U_{r_i},U_{r(\s(i))}}^{\beta}=P^{\beta}_{\eta,N,\Ucal},\qquad \eta\in\Mcal_1^{\ssup{{\rm s},N}}(\Sigma^2),\sigma\in \Sym_N({R},\eta), r_1,\dots,r_N\in\Sigma.
\end{equation} 
Note that the measure in \eqref{low4} does not depend on $\sigma$, as long as $\sigma\in \Sym_N({R},\eta)$.
Furthermore, note that $\overline \eta$ is the empirical measure of the configuration  $R=(r_1,\dots,r_N)$, and therefore
\begin{equation}\label{low5}
\prod_{i=1}^N m(r_i)=\prod_{r\in \Sigma}m(r)^{N\overline{\eta}(r)}.
\end{equation} 

On the right hand side of \eqref{low2}, we insert a sum on $\eta\in\Mcal_1^{\ssup{{\rm s},N}}(\Sigma^2)$ and restrict the sum on $\sigma$ to $\sigma\in \Sym_N({R},\eta)$. Substituting \eqref{low4} and \eqref{low5}, we arrive at 
\begin{eqnarray}\label{estlowercounting}
\begin{aligned}
\P^{\ssup{\rm sym}}_{\m,N}&(L_N\in G)\ge \sum_{\eta\in\Mcal_1^{\ssup{{\rm s},N}}(\Sigma^2)}\Big(\prod_{r\in \Sigma}m(r)^{N\overline{\eta}(r)}\Big)
P^{\beta}_{\eta,N,\Ucal}(L_N\in G_{\delta})\frac 1{N!} \sum_{R\in\Sigma^N}{\sharp\Sym_N(R,\eta)}.
\end{aligned}
\end{eqnarray}
Now we compute the counting term $ \sum_{R\in\Sigma^N}{\sharp \Sym_N(R,\eta)}$. For $R=(r_1,\dots,r_N)$ and $\s\in\Sym_N$, we write $R_\s=(r(\sigma(1)),\dots,r(\s(N)))$. Let $L(R)(r)=\frac 1{N}\sharp\{i\in\{1,\dots,N\}\colon r_i=r\}$ denote the empirical measure of the configuration $R$. In the following, we also sum over all configurations  $\psi=R_\s$ with empirical measure $L(\psi)$. Then, for any $\eta\in \Mcal_1^{\ssup{{\rm s},N}}(\Sigma^2) $, we compute
\begin{eqnarray}\label{counting}
\begin{aligned} 
\sum_{R\in\Sigma^N}{\sharp \Sym_N(R,\eta)}&=\sum_{\sigma\in\Sym_N}\sum_{R\in\Sigma^N}\sum_{{
\psi}\in\Sigma^N} 
\1\{R_{\s}={\psi}\}\1\{\forall r,s\colon \#\{i\colon r_i=r,\psi_i=s\}=N\eta(r,s)\}\\ &=\sum_{\heap{R\in\Sigma^N\colon}{\overline{\eta}=L(R)}}\sum_{\heap{{ \psi}\in \Sigma^N\colon}{\overline \eta=L(\psi)}}\1\{\forall r,s\colon \#\{i\colon r_i=r,\psi_i=s\}=N\eta(r,s)\}\sum_{\s\in\Sym_N}\1\{R_{\s}={ \psi}\}\\ 
&=\frac{N!}{\prod_{r\in\Sigma} (N\overline 
\eta(r))!}\sum_{{ \psi}\in\Sigma^N}\1\{\forall r,s\colon 
\#\{i\colon r_i=r,\psi_i=s\}=N\eta(r,s)\}\prod_{r\in\Sigma} (N\overline \eta(r))!.
\end{aligned} 
\end{eqnarray}
The last equality is explained as follows. First, it is easy to see that, for fixed $R,\psi$ having empirical measures equal to $\overline \eta$, there are precisely $\prod_{r\in\Sigma} (N\overline \eta(r))!$ permutations of the coefficients of $R$ which are indistinguishable from $\psi$. Second, observe that the term 
$$
\sum_{\heap{{ \psi}\in \Sigma^N\colon}{\overline \eta=L(\psi)}}\1\{\forall r,s\colon \#\{i\colon r_i=r,\psi_i=s\}=N\eta(r,s)\}
$$
does not depend on $R$ as long as $L(R)=\overline \eta$. It is elementary that the number of configurations $R$, whose empirical measure is equal to $\overline \eta$, is equal to $ N!/\prod_{r\in\Sigma}(N\overline{\eta}(r))! $. 

Now the remaining counting factor may be evaluated using
\begin{equation}\label{counting2}
\sum_{{ \psi}\in\Sigma^N}\1\{\forall r,s\colon 
\#\{i\colon r_i=r,\psi_i=s\}=N\eta(r,s)\}=\frac{\prod_{r\in\Sigma} (N\overline \eta(r))!}{\prod_{r,s\in\Sigma} (N \eta(r,s))!},
\end{equation}
as is e.g.~seen via a well-known formula for the number of Euler trails in a complete graph (cf.~\cite{A01} and references therein), but also follows from elementary combinatorial considerations. Thus we get from (\ref{estlowercounting})
\begin{eqnarray}\label{summarisingfast}
\begin{aligned}
\P^{\ssup{\rm sym}}_{\m,N}&(L_N\in G)\ge \sum_{\eta\in\Mcal_1^{\ssup{{\rm s},N}}(\Sigma^2)}
\Big(\prod_{r\in\Sigma}m(r)^{N\overline{\eta}(r)}\Big) \frac{\prod_{r\in\Sigma} (N\overline \eta(r))!}{\prod_{r,s\in\Sigma} (N \eta(r,s))!} P^{\beta}_{\eta,N,\Ucal}(L_N\in G_{\delta}).
\end{aligned}
\end{eqnarray} 
Using Stirling's formula, we know that there is an absolute constant $C\in(0,\infty)$ such that
\begin{equation}\label{StirlingC}
1\le\frac{N!}{(N/{\rm e})^N\sqrt{2\pi N}}\le \sqrt{\frac{C}{2\pi}}\qquad\mbox{ for any }\; N\in\N.
\end{equation}
Hence one sees that, for any $\eta\in\Mcal_1^{\ssup{{\rm s},N}}(\Sigma^2)$, 
\begin{equation}\label{Stirling}
\Big(\prod_{r\in\Sigma}m(r)^{N\overline{\eta}(r)}\Big) \frac{\prod_{r\in\Sigma} (N\overline \eta(r))!}{\prod_{r,s\in\Sigma} (N \eta(r,s))!} \geq  (CN)^{-\frac 12(\sharp\Sigma)^2} {\rm e}^{-N H(\eta|\overline{\eta}\otimes m)}.
\end{equation}
Here $H$ is the entropy made explicit in \eqref{entropy}. Substituting \eqref{Stirling} in \eqref{summarisingfast}, we arrive at the assertion.
\qed
\end{proofsect}

Now we let the partition $ \Ucal=\Ucal_N=\{U_r^{\ssup{N}}\colon r\in\Sigma_N\}$ depend on $ N $ such that the fineness $f_{\Ucal_N}$ vanishes. We also write $ m_N(r)=\frac{\m(U_r^{\ssup{N}})}{\m(\Lambda)} $ for $r\in\Sigma_N,$ and $ H_N $ for the relative entropy of pair measures on $\Sigma_N^2$.

The proof of Theorem \ref{mainnew}(i) directly follows from a combination of Proposition~\ref{firstproplower} and the following.

\begin{prop}\label{proplowerpairmeasure}
There is a sequence $ (\Ucal_N)_{N\in\N} $ of partitions $ \Ucal_N=\{U_r^{\ssup{N}}\colon r\in \Sigma_N\} $ of $\Lambda $ satisfying $ \delta_N=f_{\Ucal_N}\to 0 $ and $ \sharp\Sigma_N=o(N^{1/4})$ such that
\begin{equation}\label{prop1}
\begin{aligned}
\liminf_{N\to\infty}\frac{1}{N}&\log\Bigl(\sum_{\eta\in\Mcal_1^{\ssup{{\rm s},N}}(\Sigma_N^2)}{\rm e}^{-NH_N(\eta|\overline{\eta}\otimes m_N)}P^{\beta}_{\eta,N,\Ucal_N}(L_N\in G_{\delta_N})\Bigr)
&\ge -\inf_{\mu\in G}I_{\m}^{\ssup{\rm sym}}(\mu).
\end{aligned}
\end{equation}
\end{prop}

\begin{proofsect}{Proof} Recall \eqref{Entropydef}, \eqref{Isymdef} and \eqref{Jqdef}. It suffices to construct, for any $ \mu\in G $ and $ q\in\Mcal_1^{\ssup{{\rm s}}}((\R^d)^2) $, some $\mu_N\in\Mcal_1^{\ssup{{\rm s},N}}(\Sigma_N^2)$ such that
\begin{equation}\label{goallowbound}
\limsup_{N\to\infty}\Big(H_N(\eta_N|\overline{\eta_N}\otimes m_N)-\frac 1N\log P^{\beta}_{\eta_N,N,\Ucal_N}(L_N\in G_{\delta_N})\Big)
\leq H(q|\overline{q}\otimes \m)+I^{\ssup{q}}(\mu).
\end{equation}
Fix  $ \mu\in G $ and $ q\in\Mcal_1^{\ssup{{\rm s}}}((\R^d)^2) $. We may assume that $ H(q|\overline{q}\otimes m)<\infty $. This implies that $ \supp\,(q)\subset \Lambda\times \Lambda $, because $ H(q|\overline{q}\otimes m) =H(q|\overline{q}\otimes\overline{q})+H(\overline{q}|\m) $ and the support of $ \m $ is contained in  $ \Lambda $. 
We choose a sequence of partitions $ \Ucal_{N}=\{U_r^{\ssup{N}}\colon r\in\Sigma_N\} $ of $ \Lambda $ such that $ \sharp
\Sigma_N=o(N^{1/4})$ and $\delta_N=f_{\Ucal_N}\to 0 $ as $ N\to\infty $. We write $ U_r=U_r^{\ssup{N}} $ for any $ r\in\Sigma_N $ in the following. 

First consider $ \eta\in \Mcal_1^{\ssup{\rm s}}(\Sigma_N^2) $ defined by $ \eta(r,s)=q(U_r\times U_s), r,s\in\Sigma_N$. The sequence of probability measures $\widetilde q_N\in\Mcal_1^{\ssup{\rm s}}((\R^d)^2) $ having Lebesgue density 
$$
\frac{\widetilde q_N(\d x,\d y)}{\d x\,\d y}=\sum_{r,s\in\Sigma_N}\frac{\eta(r,s)}{|U_r\times U_s|}\1_{U_r\times U_s}(x,y), 
$$
is easily seen to converge weakly towards $ q $ as $ N\to\infty $. The main technical task consists in finding a measure $\eta_N$ in $\Mcal_1^{\ssup{{\rm s},N}}(\Sigma_N^2)$ that approximates $\eta$ well enough:

\begin{lemma}\label{lem-etaconstr} Let $ (\Ucal_N)_{N\in\N} $ be a sequence of partitions $ \Ucal_N=\{U_r^{\ssup{N}}\colon r\in \Sigma_N\} $ of $\Lambda $ satisfying $ \delta_N=f_{\Ucal_N}\to 0 $ and $ \sharp\Sigma_N=o(N^{1/4})$. Then, for any $N$ sufficiently large and for any $\eta\in \Mcal_1^{\ssup{\rm s}}(\Sigma_N^2) $, there is $\eta_N\in\Mcal_1^{\ssup{{\rm s},N}}(\Sigma_N^2)$ such that 
\begin{equation}\label{etaappr}
\max_{r,s\in\Sigma_N}\big|\eta(r,s)-\eta_N(r,s)\big|\leq 2\frac {(\sharp \Sigma_N)^2}N.
\end{equation}
\end{lemma}

\begin{proofsect}{Proof}
Due to the assumption that $ \sharp\Sigma_N=o(N^{1/4}) $, we may assume that there is a $ (r_0,s_0)\in\Sigma_N^2 $ with
\begin{equation}\label{assumponeta}
\eta(r_0,s_0)\ge 2\frac{(\sharp\Sigma_N)^2}{N}\quad\mbox{ for all } N\in\N,
\end{equation}
because, if otherwise all entries are strictly smaller than $ 2(\sharp\Sigma_N)^2/N $, then $\eta$ would have, for all large $N$, total mass $ \sum_{r,s\in\Sigma_N} \eta(r,s)<2(\sharp\Sigma_N)^2(\sharp\Sigma_N)^2/N<1 $. Without loss of generality, we assume that $ r_0\not= s_0 $. The case $ r_0=s_0 $ is in fact easier and follows analogously. 
We denote by $ \floor{x} $ the largest integer smaller or equal to $x\in\R_+ $. We define $ \eta_N\colon\Sigma_N\times\Sigma_N\to\R $ by
\begin{eqnarray}
\eta_N(r,s)&=&\frac{\floor{N\eta(r,s)}}{N}\quad\mbox{ for }\; r\in\Sigma_N\setminus\{r_0\}, s\in\Sigma_N\setminus\{s_0\};\label{firstentries}\\
\eta_N(r,s_0)&=&\frac{\floor{N\overline{\eta}(r)}}{N}-\sum_{s\in\Sigma_N\setminus\{s_0\}}\eta_N(r,s)\quad\mbox{ for }\; r\in\Sigma_N\setminus\{r_0\};\label{etaN2}\\
\eta_N(r_0,s)&=&\frac{\floor{N\overline{\eta}(s)}}{N}-\sum_{r\in\Sigma_N\setminus\{r_0\}}\eta_N(r,s)\quad\mbox{ for }\; s\in\Sigma_N\setminus\{r_0,s_0\};\label{etaN3}\\
\eta_N(r_0,r_0)&=&1-\sum_{r\in\Sigma_N\setminus\{r_0\}}\frac{\floor{N\overline{\eta}(r)}}{N}-\sum_{r\in\Sigma_N\setminus\{r_0\}}\eta_N(r,r_0);\label{etaN4}\\
\eta_N(r_0,s_0)&=&1-\sum_{\heap{(r,s)\in\Sigma_N^2\colon}{(r,s)\not= (r_0,s_0)}}\eta_N(r,s)\label{etaN5}.
\end{eqnarray}

Obviously, $ \eta_N(r,s)\in \frac 1N \N_0$ for any $r,s\in\Sigma_N$. Furthermore, by \eqref{etaN5}, they sum up to one. It remains to show that $\eta_N(r,s)\geq 0$ for any $r,s\in\Sigma_N$, that $\eta_N$ satisfies the marginal property, and that \eqref{etaappr} holds.

From \eqref{firstentries} it is clear that $ 0\leq \eta_N(r,s)\leq\eta(r,s) $ for $ r\in\Sigma_N\setminus\{r_0\}$ and $ s\in\Sigma_N\setminus\{s_0\} $. Using the estimate $ \floor{x+y}\ge \floor{x}+\floor{y} $ for any $x,y\in\R_+,$ we see from \eqref{etaN2} and \eqref{etaN3} that $ \eta_N(r,s_0)\geq 0 $ for $ r\in\Sigma_N\setminus\{r_0\} $ as well as $ \eta_N(r_0,s)\geq 0 $ for $ s\in\Sigma_N\setminus\{r_0,s_0\} $. Using $ \floor{x}\le x $ we estimate
$$
\eta_N(r_0,r_0)\ge 1-\sum_{r\in\Sigma_N\setminus\{r_0\}}\overline{\eta}(r)-\sum_{r\in\Sigma_N\setminus\{r_0\}}\eta(r,r_0)=\overline{\eta}(r_0)-\sum_{r\in\Sigma_N\setminus\{r_0\}}\eta(r,r_0)\ge 0.
$$
Using now the estimate $ x-1/N\le \floor{Nx}/N\le x $ for $ x\in\R_+ $, we see from \eqref{etaN2} that $\eta_N(r,s_0)\leq \eta(r,s)+(\sharp\Sigma_N-1)/N$ for any $ r\in\Sigma_N\setminus\{r_0\}$, and  we see from \eqref{etaN3} that
$\eta_N(r_0,s)\leq \eta(r_0,s)+(\sharp\Sigma_N-1)/N$ for any $ s\in\Sigma_N\setminus\{r_0,s_0\}$. In the same way, we see from \eqref{etaN4} that
$$
\eta_N(r_0,r_0)\leq 1-\sum_{r\in\Sigma_N\setminus\{r_0\}}\overline\eta(r)+2\frac{\sharp\Sigma_N-1}N -\sum_{r\in\Sigma_N\setminus\{r_0\}}\eta(r,r_0)=\eta(r_0,r_0)+2\frac{\sharp\Sigma_N-1}N.
$$
Using all the preceding estimates, we see that
$$
\sum_{\heap{(r,s)\in\Sigma_N^2\colon}{(r,s)\not= (r_0,s_0)}}\eta_N(r,s)\leq \sum_{\heap{(r,s)\in\Sigma_N^2\colon}{(r,s)\not= (r_0,s_0)}}\eta(r,s)+2\frac{(\sharp\Sigma_N-1)^2}N+2\frac{\sharp\Sigma_N-1}N,
$$
and \eqref{etaN4} implies that $\eta_N(r_0,s_0)\geq \eta(r_0,s_0)-2(\sharp\Sigma)^2/N$, which is nonnegative by \eqref{assumponeta}. Hence, we have shown that $\eta_N$ is a probability measure on $\Sigma_N\times\Sigma_N$. From the preceding, it is also clear that \eqref{etaappr} holds.

It remains to show the marginal property of $\eta_N$. The first marginals, i.e., the sum over the right entries, are identified as
\begin{equation}\label{firstmarginals}
\sum_{s\in\Sigma_N}\eta_N(r,s)=\frac{\floor{N\overline{\eta}(r)}}{N}\quad\mbox{ for }\; r\in\Sigma_N\setminus\{r_0\}\;\mbox{ and }\; \sum_{s\in\Sigma_N}\eta_N(r_0,s)=1-\sum_{r\in\Sigma_N\setminus\{r_0\}}\frac{\floor{N\overline{\eta}(r)}}{N}.
\end{equation}
We check that they coincide with the second marginals, i.e., the sums over the left entries. For $ r\in\Sigma_N\setminus\{r_0,s_0\} $ we get from \eqref{firstentries} and \eqref{etaN2} that $ \sum_{s\in\Sigma_N}\eta^{\ssup{N}}(s,r)=\frac{\floor{N\overline{\eta}(r)}}{N} $; hence the marginals coincide for   $ r\in\Sigma_N\setminus\{r_0,s_0\} $. Using \eqref{etaN4} we see that $ \sum_{r\in\Sigma_N}\eta_N(r,r_0)=\eta_N(r_0,r_0)+\sum_{r\in\Sigma_N\setminus\{r_0\}}\eta_N(r,r_0)=1-\sum_{r\in\Sigma_N\setminus\{r_0\}}\frac{\floor{N\overline{\eta}(r)}}{N} $, and hence the marginals coincide also in $ r_0 $. Since all marginals of $\eta_N$ are probability measures on $\Sigma_N$, the two marginals coincide also in $s_0$. This shows that $\eta_N\in\Mcal_1^{\ssup{{\rm s},N}}(\Sigma_N^2) $.
\qed\end{proofsect}

Let $\eta_N$ be as in Lemma~\ref{lem-etaconstr} for $ \eta$ defined by $ \eta(r,s)=q(U_r\times U_s)$ as above. Consider the probability measures $ q_N\in\Mcal_1^{\ssup{\rm s}}(\R^d\times\R^d) $, having Lebesgue density 
$$
\frac{q_N(\d x,\d y)}{\d x\,\d y}=\sum_{r,s\in\Sigma_N}\frac{\eta_N(r,s)}{|U_r\times U_s|}\1_{U_r\times U_s}(x,y).
$$ 
From \eqref{etaappr} and the convergence of $\widetilde q_N$ towards $q$ we have that also $q_N$
converges  weakly towards $ q$. To see this, note that for any $ g\in\Ccal_{\rm b}(\R^d\times\R^d) $,
\begin{equation}\label{esttildeetaN}
\begin{aligned}
\big|\langle g,q\rangle-\langle g,q_N\rangle \big|&\leq\big|\langle g,q\rangle-\langle g,\widetilde q_N\rangle \big|+ \sum_{r,s\in\Sigma_N}\int_{U_s}\int_{U_r}\frac{|\eta(r,s)-\eta_N(r,s)|}{|U_r\times U_s|}|g(x,y)|\,\d x\d y\\
&\le o(1) +2 ||g||_{\infty}\sum_{r,s\in\Sigma_N}\frac{(\sharp\Sigma_N)^2}{N}\le o(1) +2||g||_{\infty}\frac{(\sharp\Sigma_N)^4}{N}=o(1),\qquad\mbox{as }N\to\infty.
\end{aligned}
\end{equation}

Observe that the marginals $ \overline{q}_N $ of $ {q}_N $ have Lebesgue density $x\mapsto \sum_{r\in\Sigma_N}\overline\eta_N(r)\1_{U_r}(x)/|U_r|$; in particular $\overline{q}_N(U_r)=\overline\eta_N(r)$.
Note further that the relative entropy can be written as
\begin{equation}\label{entropysum}
H(q_N|\overline{q}_N\otimes\m)=H(\overline{q}_N|\m)+H(q_N|\overline{q}_N\otimes\overline{q}_N).
\end{equation}
Jensen's inequality, applied for the function $\varphi(z)=z\log z$, gives for the first entropy on the right of \eqref{entropysum},
\begin{equation}
\begin{aligned}
H(\overline{q}_N|\m)&=\int_{\R^d}\d\m\,\frac{\d\overline{q}_N}{\d\m}\log\frac{\d\overline{q}_N}{\d\m}
=\sum_{r\in \Sigma_N}\m(U_r)\int_{U_r}\frac{\d\m}{\m(U_r)}\,\varphi\Big(\frac{\d\overline{q}_N}{\d\m}\Big)\\
&\ge\sum_{r\in\Sigma_N}\m(U_r)\varphi\Big(\int_{U_r}\frac{\d\m}{\m(U_r)}\,\frac{\d\overline{q}_N}{\d\m}\Big)=\sum_{r\in\Sigma_N} \overline{q}_N(U_r)\log\frac{\overline{q}_N(U_r)}{\m(U_r)}\\
&=\sum_{r\in\Sigma_N}\overline{\eta}_N(r)\log\frac{\overline{\eta}_N(r)}{m_N(r)}=H_N(\overline{\eta}_N|m_N).
\end{aligned}
\end{equation}
In the same way one shows that $H(q_N|\overline{q}_N\otimes\overline{q}_N)\geq H_N(\eta_N|\overline{\eta}_N\otimes\overline{\eta}_N)$,
resulting in
$$
H(q_N|\overline{q}_N\otimes\m)\ge H_N(\eta_N|\overline{\eta}_N\otimes m_N).
$$
By \cite[Prop.~15.6]{Geo88}, we have $\lim_{N\to\infty}H(q_N|\overline{q_N}\otimes\m)=H(q|\overline{q}\otimes\m)$. Hence, we have shown that
$$
\limsup_{N\to\infty}H_N(\eta_N|\overline{\eta_N}\otimes m_N)\leq H(q|\overline{q}\otimes\m).
$$
That is, we have shown the first half of \eqref{goallowbound}. 

In our final step, we are going to use the G\"artner-Ellis Theorem to deduce that 
$$ 
\liminf_{N\to\infty}\frac{1}{N}\log P^{\beta}_{\eta_N,N,\Ucal_N}(L_N\in G_{\delta}) \ge - I^{\ssup{q}}(\mu). 
$$ 
For doing this, we first introduce, for any $ \Phi \in\Ccal_{\rm b}(\Ccal) $,
\begin{equation}
\begin{aligned}
\Lcal_N(\Phi )&:=\log \E^{\beta}_{\eta_N,N,\Ucal_N}\bigl[{\rm e}^{N\langle \Phi ,L_N\rangle}\bigr]
=\log\Bigl(\prod_{r,s\in\Sigma_N}\E^{\beta}_{U_r,U_s}\big[{\rm e}^{\Phi (B)}\big]^{N\eta_N(r,s)}\Bigr)\\
&=N\sum_{r,s\in\Sigma_N}\eta_N(r,s)\log \E^{\beta}_{U_r,U_s}\big[{\rm e}^{\Phi (B)}\big]\\
&=N\int_{\R^d}\int_{\R^d}q_N(\d x,\d y)\log\E^{\beta}_{U_{r_N(x)},U_{r_N(y)}}\big[{\rm e}^{\Phi(B)}\big],
\end{aligned}
\end{equation} 
where $ r_N(x)\in\Sigma_N $ is defined by $ x\in U_{r_N(x)}$. From the proof of Lemma~\ref{estimationsetendpoints} it is easily seen that 
$$ 
\lim_{N\to\infty} \E^{\beta}_{U_{r_N(x)},U_{r_N(y)}}\big[{\rm e}^{\Phi (B)}\big]=\E^{\beta}_{x,y}\big[{\rm e}^{\Phi (B)}\big], 
$$ 
uniformly in $ x,y\in \Lambda $. Recall that $ q_N\to q $ as $ N\to\infty $ weakly. Hence, the limit $ \Lcal(\Phi )=\lim_{N\to\infty}\frac 1N\Lcal_N(\Phi ) $ exists, and 
$$
\Lcal(\Phi )=\int_{\R^d}\int_{\R^d}q(\d x,\d y)\,\log\E^{\beta}_{x,y}\big[{\rm e}^ {\Phi (B)}\big].
$$ 
Since it is easily seen that $ \Lcal $ is lower semi continuous and  G\^ateaux differentiable,  and by the exponential tightness of the family $ (P_{\eta_N,N,\Ucal_N}^{\beta})_{N\in\N} $ (see Lemma~\ref{tightnessPeta}), \cite[4.5.27]{DZ98} implies that 
\begin{equation}
\liminf_{N\to\infty}\frac{1}{N}\log P^{\beta}_{\eta_N,N,\Ucal_N}(L_N\in G_{\delta})\ge - I^{\ssup{q}}(\mu).
\end{equation}
This shows the second half of \eqref{goallowbound} and ends the proof.
\qed
\end{proofsect}

\subsection{Proof of Theorem~\ref{mainnew}(ii).}\label{sec-uppproof} 

\noindent Our proof of the upper bound uses the same machinery as the proof of the lower bound.
Recall that we assume that $ \m $ is a probability measure on $ \R^d $, not necessarily having compact support. Also recall the notation from the beginning of Section~\ref{sec-lowproof}, in particular, the partition  $ \Ucal=\{U_r\colon r\in\Sigma\} $ of a given set $ \Lambda\subset\R^d $ and \eqref{pairmeasures}--\eqref{entropy}. 

In the following, we will have to work with probability measures on $\Sigma\times \Sigma$ that satisfy the marginal property only approximatively. For $\eps\in(0,1)$ and $n\in\N$ introduce the set

$$
\Mcal_1^{\ssup{\eps}}(\Sigma^2)=\Big\{\eta\in \Mcal_1(\Sigma^2)\colon \d(\eta^{\ssup 1},\eta^{\ssup 2})\leq 2\eps\Big\},
$$
where  $\d$ is some metric on $\Mcal_1(\Sigma^2)$ that induces the weak topology, and $ \eta^{\ssup{1}} $ and $ \eta^{\ssup 2}  $ are the two marginal measures of $ \eta $. By $\Mcal_1^{\ssup{\eps,n}}(\Sigma^2)$ we denote the set $\Mcal_1^{\ssup{\eps}}(\Sigma^2)\cap\frac 1n\N_0^{\Sigma^2}$.

Our first main step is the following.

\begin{prop}[Combinatorics]\label{propcombupper}Fix a closed set  $ F\subset\Mcal_1(\Ccal) $ and a compact set $ \Lambda\subset\R^d $. Then, for any $ \delta>0 $, any partition $ \Ucal $ of $ \Lambda $ having fineness $ f_{\Ucal}\le\delta $  and for any $ \eps >0 $ and $ N\in\N $,
\begin{equation}\label{countingupper5}
\begin{aligned}
\P_{\m,N}^{\ssup{\rm sym}}(L_N\in F)&\le 2^N \m(\Lambda^{\rm c})^{\eps N}+(CN)^{\frac 12(\sharp\Sigma)^2}{\rm e}^{NC_{\eps}}\sum_{(1-2\eps)N<n\leq N}\\
&\qquad\times\sum_{\eta\in \Mcal^{\ssup{\eps,n}}(\Sigma^2)}{\rm e}^{-n H(\eta| m\otimes \eta^{\ssup 2})} P^{\beta}_{\eta,n,\Ucal}(L_{n}\in F^{2\eps+\delta}),
\end{aligned}
\end{equation}
where $ C>0 $ is given in \eqref{StirlingC}, and $ C_{\eps}>0 $ vanishes as $ \eps\downarrow 0$. 
\end{prop}

\begin{proofsect}{Proof}
Consider \eqref{Psymdef}. We split each of the $ N $ integrations over the starting points of the Brownian bridges into an integration over $ \Lambda $ and over the complement $ \Lambda^{\rm c} $. Thus we can write
\begin{equation}\label{splittingint1}
\int_{(\R^d)^N}=\sum_{a\in\{1,{\rm c}\}^N}\int_{\Lambda^{a_1}}\cdots\int_{\Lambda^{a_N}},
\end{equation}
where we used the notation $\Lambda^1=\Lambda$.
The sum on $a$ is split into the two sums where more than $ \eps N $ integrals are on $\Lambda$ and the remainder:
\begin{equation}\label{splittingint2}
\int_{(\R^d)^N}=\sum_{\heap{a\in\{1,{\rm c}\}^N\colon}{\sharp\{i\colon a_i={\rm c}\}\ge \eps N}}\int_{\Lambda^{a_1}}\cdots\int_{\Lambda^{a_N}} +\sum_{\heap{a\in\{1,{\rm c}\}^N\colon}{ \sharp\{i\colon a_i=1\}> (1-\eps) N}}\int_{\Lambda^{a_1}}\cdots\int_{\Lambda^{a_N}}.
\end{equation} 
Using this in \eqref{Psymdef}, we write $\P_{\m,N}^{\ssup{\rm sym}}=\P_{\m,N}^{\ssup{\rm sym},I}+\P_{\m,N}^{\ssup{\rm sym},II}$, with obvious notation. It is clear that
\begin{equation}\label{PIupperesti}
\P_{\m,N}^{\ssup{\rm sym},I}(L_N\in F)\leq 2^N \m(\Lambda^{\rm c})^{\eps N}.
\end{equation}
This is the first term of the right hand side of \eqref{countingupper5}, and now we show that the second part is an estimate for $\P_{N,\beta}^{\ssup{\rm sym},II}(L_N\in F)$.
For doing this, we distinguish all the sets $I$ of indices $i$ such that $a_i=1$:
\begin{equation}\label{upperest1}
\begin{aligned}
\P_{\m,N}^{\ssup{\rm sym},II}(L_N\in F)
&= \frac{1}{N!}\sum_{\s\in\Sym_N} \sum_{\heap{I\subset\{1,\ldots,N\}}{\sharp I>(1-\eps)N}}\sum_{\heap{a\in\{1,{\rm c}\}^N}{I=\{i\colon a_i=1\}}}\\
&\qquad\qquad\qquad\int_{\Lambda^{a_1}}\cdots\int_{\Lambda^{a_N}}\m(\d x_1)\cdots\m(\d x_N)\Bigl(\bigotimes_{i=1}^N\P_{x_i,x_{\s(i)}}^{\beta}\Bigr)(L_N\in F).
\end{aligned}
\end{equation}
In the product over the Brownian bridges we only want to consider those Brownian bridges whose initial and terminal points are in $ \Lambda $. Given $ \s\in\Sym_N $ and $ I\subset\{1,\ldots,N\} $, we consider the subset $ I_{\s}=\{i\in \{1,\dots,N\}\colon \s(i)\in I\}=\s^{-1}(I) $. We want to replace the measure $ \bigotimes_{i=1}^N\P_{x_i,x_{\s(i)}}^{\beta} $ by the measure $ \bigotimes_{i\in I_{\s}\cap I}\P_{x_i,x_{\s(i)}}^{\beta} $, i.e., we want to forget about all the motions whose initial point $x_i$ or whose terminal point $x_{\s(i)}$ is not in $\Lambda$. We do this by replacing the empirical path measure $ L_N $ by the empirical path measure $L_{I_\sigma\cap I}$, where
$$
L_{J}=\frac{1}{\sharp J}\sum_{i\in J} \delta_{B^{\ssup{i}}},\qquad J\subset\{1,\dots,N\}.
$$
Recall that we work with sets $I\subset\{1,\dots,N\}$ satisfying $\sharp I>(1-\eps)N$ and have therefore  $ \d (L_N,L_{I_{\s}\cap I})< 2\eps $, since $\sharp(I_{\s}\cap I)>(1-2\eps)N$. Hence, if $F^{\delta}=\{\mu\in \Mcal_1(\Ccal)\colon\dist(\mu,F)\leq \delta\}$ denotes the closed $\delta$-neighbourhood of $F$, then
\begin{equation}
\sharp I>(1-\eps)N\qquad\Longrightarrow\qquad\Bigl(\bigotimes_{i=1}^N\P_{x_i,x_{\s(i)}}^{\beta}\Bigr)(L_N\in F)\le \Bigl( \bigotimes_{i\in I_{\s}\cap I}\P_{x_i,x_{\s(i)}}^{\beta}\Bigr)(L_{I_{\s}\cap I}\in F^{2\eps}).
\end{equation}
Using this in \eqref{upperest1}, we can freely execute the $ N-\sharp I $ integrations over those $\m(\d x_j)$ with $j\notin I$ since they do not contribute anymore. These integrations may be estimated from above by one, and we are left with the $\sharp I$ integrations over those $x_i$ satisfying $i\in I$, which means that $x_i\in \Lambda$. Hence, all the remaining integration areas are equal to $\Lambda$. Note that then the sum on all $a\in\{1,{\rm c}\}^N$ satisfying $I=\{i\colon a_i=1\}$ just yields a factor of one. This gives
\begin{equation}\label{upperest2}
\begin{aligned}
&\P_{\m,N}^{\ssup{\rm sym},II}(L_N\in F)\le \frac{1}{N!}\sum_{\s\in\Sym_N}\sum_{\heap{I\subset\{1,\ldots,N\}}{\sharp I>(1-\eps)N}}\int_{\Lambda^{I}}\prod_{i\in I}\m(\d x_i)\Bigl(\bigotimes_{i\in I_{\s}\cap I}\P_{x_i,x_{\s(i)}}^{\beta}\Bigr)(L_{I_{\s}\cap I}\in F^{2\eps}).
\end{aligned}
\end{equation} 
Now we introduce a partition $\Ucal=\{U_r\colon r\in\Sigma\}$ of $\Lambda$ and split the integration over $\Lambda^I$ into a sum on integrations like in \eqref{split}:
\begin{equation}\label{splitupper}
\int_{\Lambda^{I}}\prod_{i\in I}\m(\d x_i)=\sum_{R\in\Sigma^{I}}\prod_{i\in I}\int_{U_{r(i)}}\m(\d x_i),\qquad R=(r(i))_{i\in I}.
\end{equation}
For fixed $R\in \Sigma^I$ and for multi-indices $x_i\in U_{r(i)}$ with $i\in I$, we may estimate
\begin{equation}\label{splitupper2}
\Bigl(\bigotimes_{i\in I_{\s}\cap I}\P_{x_i,x_{\s(i)}}^{\beta}\Bigr)(L_{I_{\s}\cap I}\in F^{2\eps})\le \sup_{\heap{y_i\in U_{r(i)},}{i\in I}}\Bigl(\bigotimes_{i\in I_{\s}\cap I}\P_{x_{i},y_{\s(i)}}^{\beta}\Bigr)(L_{I_{\s}\cap I}\in F^{2\eps}).
\end{equation}
The right hand side does not depend on the $x_i$ with $i\in I\setminus I_\s$. Hence, after substituting \eqref{splitupper} and \eqref{splitupper2} in \eqref{upperest2}, the integrations over $x_i\in U_{r(i)}$ with $i\in I\setminus I_\s$ may be executed freely and their contribution gives a factor of $m(r(i))$, where we recall that $m(r)=\m(U_r)$ for $r\in \Sigma$.
Now we perform the integrations over all the remaining $x_i$, i.e., over $x_i\in U_{r(i)}$ with $i\in I_{\s}\cap I$. Recall the notation in \eqref{PUydef}, to obtain, also using \eqref{splitupper} and \eqref{splitupper2},
\begin{equation}
\begin{aligned}
\int_{\Lambda^{I}}&\prod_{i\in I}\m(\d x_i)\Bigl(\bigotimes_{i\in I_{\s}\cap I}\P_{x_i,x_{\s(i)}}^{\beta}\Bigr)(L_{I_{\s}\cap I}\in F^{2\eps})\\
&\leq \sum_{R\in\Sigma^{I}}\prod_{i\in I}m(r(i))\sup_{\heap{y_i\in U_{r(i)},}{i\in I}}\Bigl(\bigotimes_{i\in I_{\s}\cap I}\P_{U_{r(i)},y_{\s(i)}}^{\beta}\Bigr)(L_{I_{\s}\cap I}\in F^{2\eps})\\
&\leq \sum_{R\in\Sigma^{I}}\prod_{i\in I}m(r(i))\Bigl(\bigotimes_{i\in I_{\s}\cap I}\P_{U_{r(i)},U_{r(\s(i))}}^{\beta}\Bigr)(L_{I_{\s}\cap I}\in F^{2\eps+\delta}).
\end{aligned}
\end{equation}
In the last step, we introduced some small $\delta>0$, assumed that the fineness $f_\Ucal=\max_{r\in \Sigma}\diam(U_r)$ is smaller than $\delta$, and used Lemma~\ref{estimationsetendpoints} (recall the notation in \eqref{PUUdef}). So far, we have deduced that
\begin{equation}
\P_{\m,N}^{\ssup{\rm sym},II}(L_N\in F)
\leq \frac{1}{N!}\sum_{\s\in\Sym_N} \sum_{\heap{I\subset\{1,\ldots,N\}}{|I|>(1-\eps)N}}\sum_{R\in\Sigma^{I}}\prod_{i\in I}m(r(i))\Bigl(\bigotimes_{i\in I_{\s}\cap I}\P_{U_{r(i)},U_{r(\s(i))}}^{\beta}\Bigr)(L_{I_{\s}\cap I}\in F^{2\eps+\delta}).
\end{equation}
Put $n=\sharp (I_{\s}\cap I)$. Observe that the  probability measure $\bigotimes_{i\in I_{\s}\cap I}\P_{U_{r(i)},U_{r(\s(i))}}^{\beta}$ does not depend on the full information contained in $\s$, but only on the frequencies of $i\in I$ such that $r(i)=r$ and $r(\s(i))=s$, for any $r,s\in \Sigma$. 
In the next step we add a sum over pair measures $ \eta$ in $\Mcal_1^{\ssup{n}}(\Sigma^2) $, the set of probability measures $\Sigma^2\to\frac 1{n}\N_0 $, and add the constraint that these frequencies are equal to $n\eta(r,s)$. Under this constraint, $\bigotimes_{i\in I_{\s}\cap I}\P_{U_{r(i)},U_{r(\s(i))}}^{\beta}$ does not depend on $\s$, but only on $\eta$, such that we may just count all the $\s$'s that satisfy the constraint. Note that these $\eta$'s do not necessarily have equal marginals. More precisely, their left marginal $\eta^{\ssup 1}$ is equal to the empirical measure of $(r_i)_{i\in I_\s\cap I}$, and its right marginal $\eta^{\ssup 2}$ is equal to the empirical measure of $(r_{\s(i)})_{i\in I_\s\cap I}$. However, since $n>(1-2\eps)N$, these $\eta$'s are elements of the set $\Mcal_1^{\ssup{\eps,n}}$ defined prior to the lemma. Thus,
\begin{equation}\label{upperest3}
\begin{aligned}
&\P_{\m,N}^{\ssup{\rm sym},II}(L_N\in F)\le 
\sum_{\heap{\widetilde I\subset I\subset\{1,\ldots,N\};}{n=\sharp \widetilde I>(1-2\eps)N}}
\sum_{\eta\in \Mcal_1^{\ssup{\eps,n}}(\Sigma^2)} \prod_{r\in\Sigma}m(r)^{n {\eta}^{\ssup 1}(r)}\Bigl(\bigotimes_{r,s\in\Sigma}\Bigl(\P_{U_r,U_s}^{\beta}\Bigr)^{n\eta(r,s)}\Bigr)(L_{n}\in F^{2\eps+\delta})\\
&\quad\times \sum_{R\in\Sigma^{I}}\sum_{\s\in\Sym_N}
\1\{\widetilde I=I_{\s}\cap I\}\frac 1{N!} \1\{\s\in\Sym_N\colon \forall\;r,s\colon\sharp\{i\in \widetilde I\colon r_i=r, r_{\s(i)}=s\}=n\eta(r,s)\}\prod_{i\in I\setminus\widetilde I}m(r_i)\\
&\le 
\sum_{(1-2\eps)N<n\leq N}
\sum_{\eta\in \Mcal_1^{\ssup{\eps,n}}(\Sigma^2)} \prod_{r\in\Sigma}m(r)^{n {\eta}^{\ssup 1}(r)}P_{\eta,n,\Ucal}^\beta (L_{n}\in F^{2\eps+\delta})\\
&\qquad\times \sum_{\heap{\widetilde I\subset I\subset\{1,\ldots,N\};}{\sharp\widetilde I=n}}
\sum_{R\in\Sigma^{I}}
\frac 1{N!} \sharp\Sym_N(R,\widetilde I,\eta)\prod_{i\in I\setminus\widetilde I}m(r_i),
\end{aligned}
\end{equation} 
where we used the notation in \eqref{low3} and introduced
$$
\Sym_N(R,\widetilde I,\eta)=\big\{\s\in\Sym_N\colon \forall\;r,s\colon\sharp\{i\in \widetilde I\colon r_i=r, r_{\s(i)}=s\}=n\eta(r,s)\}.
$$
Let us estimate the combinatorial terms in the last line of \eqref{upperest3} as follows. Fix  $(1-2\eps)N<n\leq N$, $\eta\in \Mcal^{\ssup{\eps,n}}(\Sigma^2)$ and $\widetilde I\subset I\subset\{1,\ldots,N\}$ satisfying $\sharp\widetilde I=n$. Furthermore, fix  $R=(r_i)_{i\in I}\in\Sigma^{I}$ such that $\eta^{\ssup 1}$ is equal to the empirical measure of $(r_i)_{i\in \widetilde I}$. We add a sum on $\psi\in\Sigma^{\widetilde I}$ with the constraint that $\psi=(r_{\s(i)})_{i\in \widetilde I}$. Recall that $L(\psi)\in \Mcal_1^{\ssup n}(\Sigma)$ is the empirical measure of $\psi$. This gives
\begin{equation}\label{countingupper1}
\begin{aligned}
\sharp\Sym_N(R,\widetilde I,\eta)
&=\sum_{\heap{\psi\in\Sigma^{\widetilde I}}{L(\psi)=\eta^{\ssup 2}}}\1\{\forall\, r,s\colon\sharp\{i\in \widetilde I\colon r_i=r;\psi_i=s\}=n\eta(r,s)\}\sum_{\s\in\Sym_N}\1\{\psi_i=r_{\s(i)},\,\forall\; i\in \widetilde I\}.
\end{aligned}
\end{equation}
The last term is written and estimated as follows:
\begin{equation}\label{countingupper2}
\begin{aligned}
\sum_{\s\in\Sym_N}&\1\{\psi_i=r_{\s(i)},\,\forall\; i\in \widetilde I\}
=\sum_{I'\subset I\colon \sharp I'=\sharp \widetilde I}\sharp\big\{\s\colon\widetilde I\to I'\text{ bijective}\colon \psi_i=r_{\s(i)},\,\forall\; i\in \widetilde I\big\}\\
&\qquad\times \sharp\big\{\s\colon\{1,\dots,N\}\setminus\widetilde I\to \{1,\dots,N\}\setminus I'\text{ bijective}\big\}\\
&\leq \binom{\sharp I}n\Big(\prod_{r\in\Sigma}\big(n\eta^{\ssup 2}(r)\big)!\Big)(N-n)!.
\end{aligned}
\end{equation}
 The remaining term in \eqref{countingupper1} is identified as
\begin{equation}\label{countingupper3}
\begin{aligned}
\sum_{\heap{\psi\in\Sigma^{\widetilde I}}{L(\psi)=\eta^{\ssup 2}}}\1\{\forall\, r,s\colon\sharp\{i\in \widetilde I\colon r_i=r;\psi_i=s\}=n\eta(r,s)\}
&=\frac{\prod_{r\in\Sigma}\big(n\eta^{\ssup 1}(r)\big)!}{\prod_{r,s\in\Sigma}\big(n\eta(r,s)\big)!},
\end{aligned}
\end{equation}
which is derived in the same way as \eqref{counting2} above. Note that the estimates in \eqref{countingupper2}--\eqref{countingupper3} do not depend on $R$ as long as $\eta^{\ssup 1}$ is equal to the empirical measure of $(r_i)_{i\in \widetilde I}$. The number of these $ (r_i)_{i\in I\setminus\widetilde I} $ is equal
$ n!/\prod_{r\in\Sigma}(n\eta^{\ssup 1}(r))!$. We write the sum on $ R\in\Sigma^I $ as sums on $ (r_i)_{i\in\widetilde I}\in \Sigma^{\widetilde I} $ and on $  (r_i)_{i\in I\setminus\widetilde I}\in \Sigma^{I\setminus \widetilde I} $. Taking into account the term $ \prod_{i\in I\setminus\widetilde I}m(r_i) $ in the last line of \eqref{upperest3}, the latter sum can be estimated against one. Hence, the last line of \eqref{upperest3} can be estimated as follows.
\begin{equation}\label{countingupper4}
\begin{aligned}
\sum_{\heap{\widetilde I\subset I\subset\{1,\ldots,N\};}{\sharp\widetilde I=n}}
&\sum_{R\in\Sigma^{I}}
\frac 1{N!} \sharp\Sym_N(R,\widetilde I,\eta)\prod_{i\in I\setminus\widetilde I}m(r_i)
\leq \sum_{\heap{\widetilde I\subset I\subset\{1,\ldots,N\};}{\sharp\widetilde I=n}}\sharp\big\{(r_i)_{i\in \widetilde I}\in\Sigma^{\widetilde I}\colon L\big((r_i)_{i\in \widetilde I}\big)=\eta^{\ssup 1}\big\}\\
&\qquad\times\binom{\sharp I}n\frac{\prod_{r\in\Sigma}\big(n\eta^{\ssup 1}(r)\big)!\,\prod_{r\in\Sigma}\big(n\eta^{\ssup 2}(r)\big)!}{\prod_{r,s\in\Sigma}\big(n\eta(r,s)\big)!}\frac{(N-n)!}{N!}\\
&=\binom {\sharp I}n^3 N \frac{n!}{\prod_{r\in\Sigma}\big(n\eta^{\ssup 1}(r)\big)!}\frac{\prod_{r\in\Sigma}\big(n\eta^{\ssup 1}(r)\big)!\,\prod_{r\in\Sigma}\big(n\eta^{\ssup 2}(r)\big)!}{\prod_{r,s\in\Sigma}\big(n\eta(r,s)\big)!}\frac{(N-n)!}{N!}\\
&\leq \binom Nn^2N \frac{\prod_{r\in\Sigma}\big(n\eta^{\ssup 2}(r)\big)!}{\prod_{r,s\in\Sigma}\big(n\eta(r,s)\big)!}.
\end{aligned}
\end{equation}
In \eqref{upperest3}, we substitute \eqref{countingupper4} and use Stirling's formula in a similar way as in \eqref{Stirling}, to arrive at the assertion.
\qed
\end{proofsect}

Now we use large-deviation arguments for identifying the large-deviation rate of the last line of \eqref{countingupper5}. Introduce the rate function
 \begin{equation}\label{JSigmadef}
\begin{aligned}
I_\Ucal^{\ssup\eta}(\mu)=\sup_{\Phi\in\Ccal_{\rm b}(\Ccal)}\Big(\langle \Phi,\mu\rangle-\sum_{r,s\in\Sigma}\eta(r,s)\log \E^\beta_{U_r,U_s}\big[{\rm e}^{\Phi(B)}\big]\Big).
\end{aligned}
\end{equation}

\begin{lemma}[Large deviations]\label{LDPargupper}Fix a closed set $F\subset \Mcal_1(\Ccal)$ and a compact set $ \Lambda\subset\R^d$. Then, for any $\eps,\delta>0$, and any partition $\Ucal$ of $\Lambda$ with fineness $\leq \delta$,
\begin{equation}\label{countingupper6}
\begin{aligned}
\limsup_{n \to\infty}&\frac 1n\log \Big(\sum_{\eta\in \Mcal_1^{\ssup{\eps,n}}(\Sigma^2)}{\rm e}^{-n H(\eta| m\otimes \eta^{\ssup 2})} P_{\eta,n,\Ucal}^\beta(L_{n}\in F^{2\eps+\delta})\Big)\\
&\leq  -\inf_{\mu\in F^{2\eps+\delta}} \inf_{\eta\in \Mcal_1^{\ssup{\eps}}(\Sigma^2)}\Big\{H(\eta|m\otimes \eta^{\ssup 2})+I_\Ucal^{\ssup\eta}(\mu)\Big\}.
\end{aligned}
\end{equation}
\end{lemma}

\begin{proofsect}{Proof}
From Lemma~\ref{tightnessPeta} it follows that there is a sequence of compact sets $ K_L\subset\Mcal_1(\Ccal) $ such that
\begin{equation}\label{exptightuniform}
\lim_{L\to\infty}\limsup_{n\to\infty} \frac{1}{n}\log\Bigl(\sup_{\Ucal\colon f_{\Ucal}\le\frac{1}{2}}\sup_{\eta\in\Mcal_1^{\ssup{n}}(\Sigma^2)} P_{\eta,n,\Ucal}^{\beta}(L_N\in K_L^{\rm c})\Bigr)=-\infty.
\end{equation}
Hence, it suffices to assume that $F^{2\eps+\delta}$ is a compact subset of $\Mcal_1(\Ccal)$.

We consider the logarithmic moment generating function of the distribution of $L_n$ under $P_{\eta,n,\Ucal}^\beta$,
\begin{equation}\label{logmomgenfctLN}
\begin{aligned}
\Lcal^{\ssup\eta}_{n,\Ucal}(\Phi)&=\log E_{\eta,n,\Ucal}^\beta\big[{\rm e}^{n \langle \Phi,L_n\rangle}\big]&=n\sum_{r,s\in\Sigma}\eta(r,s)\log \E^\beta_{U_r,U_s}\big[{\rm e}^{\Phi(B)}\big].
\end{aligned}
\end{equation}
Let now $\eta_n\in\Mcal_1^{\ssup{\eps,n}}(\Sigma^2)$ be maximal for $\eta\mapsto{\rm e}^{-n H(\eta| m\otimes \eta^{\ssup 2})} P_{\eta,n,\Ucal}^\beta(L_{n}\in F^{2\eps+\delta})$. Since $\Mcal_1^{\ssup{\eps}}(\Sigma^2)$ is compact, we may assume that $\lim_{n\to\infty}\eta_n=\eta$ for some $\eta_n\in\Mcal_1^{\ssup{\eps}}(\Sigma^2)$. Certainly, the limit
$$
\Lcal^{\ssup\eta}_{\Ucal}(\Phi)=\lim_{n\to\infty}\frac 1n \Lcal^{\ssup{\eta_n}}_{n,\Ucal}(\Phi)
=\sum_{r,s\in\Sigma}\eta(r,s)\log \E^\beta_{U_r,U_s}\big[{\rm e}^{\Phi(B)}\big]
$$
exists, and is lower semi continuous and G\^ateaux differentiable. Observe that $I_\Ucal^{\ssup\eta}$ is the Fenchel-Legendre transform of $\Lambda^{\ssup\eta}_{\Ucal}$. Now the G\"artner-Ellis theorem yields that 
$$
\limsup_{n\to\infty}\frac 1n \log P_{\eta_n,n,\Ucal}^\beta(L_{n}\in F^{2\eps+\delta})\leq -\inf_{\mu\in F^{2\eps+\delta}}I_\Ucal^{\ssup\eta}(\mu).
$$
Since the cardinality of $\Mcal_1^{\ssup{\eps,n}}(\Sigma^2)$ is polynomial in $n$, and by continuity of $\eta \mapsto H(\eta| m\otimes \eta^{\ssup 2})$, the assertion follows.
\qed
\end{proofsect}

Substituting Lemma~\ref{LDPargupper} on the right hand side of \eqref{countingupper5} we obtain that for any
$ \eps >0,\delta>0, $ and any compact set $ \Lambda\subset\R^d $ and any partition $ \Ucal $ of $ \Lambda $ having fineness smaller than $ \delta $
\begin{equation}\label{countingupper7}
\begin{aligned}
\limsup_{N\to\infty}\frac 1N\log&\P_{\m,N}^{\ssup{\rm sym}}(L_N\in F)\le -\min\Big\{-\log 2-\eps\log \m(\Lambda^{\rm c}),\\& C_{\eps}+\inf_{\mu\in F^{2\eps+\delta}} \inf_{\eta\in \Mcal_1^{\ssup{\eps}}(\Sigma^2)}\big\{H(\eta|m_\Ucal\otimes \eta^{\ssup 2})+I_\Ucal^{\ssup\eta}(\mu)\big\}\Big\}.
\end{aligned}
\end{equation}

In order to finish the proof of Theorem~\ref{mainnew}(ii), let $\delta\downarrow 0$ on the right hand side of \eqref{countingupper7}, replace $\eps$ and $\Lambda$ by sequences $\eps_N\downarrow 0$ and $\Lambda_N\uparrow \R^d$ such that $\eps_N\log\m(\Lambda_N^{\rm c})\to-\infty$, and use the following lemma.

\begin{lemma}\label{finenessVP} Fix a closed set $F\subset\Mcal_1(\Ccal)$. Then, for any sequence $ (\eps_N)_{N\in\N} $ in $ (0,1/2] $ satisfying $ \eps_N\to 0 $ as $ N \to\infty $ and any sequence $(\Lambda_N)_{N\in\N} $ of compact sets $ \Lambda_N\subset\R^d $ satisfying $ \Lambda_N\uparrow \R^d $ as $ N \to\infty $, 
\begin{equation}\label{countingupper8}
\begin{aligned}
\liminf_{N\to\infty}\liminf_{\delta\downarrow 0}&\inf_{\mu\in F^{2\eps_N+\delta}} \inf_{\eta\in \Mcal_1^{\ssup{\eps_N}}(\Sigma^2)}\Big\{H(\eta|m_{\Ucal_N}\otimes \eta^{\ssup 2})+I_{\Ucal_N}^{\ssup\eta}(\mu)\Big\}\\
&\geq\inf_{\mu\in F}\inf_{q\in \Mcal_1^{\ssup{{\rm s}}}(\R^d\times\R^d)}\Big\{H(q|\overline q\otimes \m)+I^{\ssup q}(\mu)\Big\},
\end{aligned}
\end{equation}
where $ \Ucal_N $ is a partition of $ \Lambda_N $ with fineness smaller or equal to $ \delta $.
\end{lemma}

\begin{proofsect}{Proof} 
Let us first roughly explain the nature of the argument. We pick approximating sequences of $\eta$'s and $\mu$'s and employ a compactness argument in order to extract a convergent subsequence. This easily finishes the proof by lower semi continuity. The compactness argument relies on the compactness of the level sets of the entropy term (which is well-known) and on that of the $I$-term, which we derive from exponential tightness of certain probability measures whose large deviation principle is governed by the $I$-term (this is in the spirit of the proof of \cite[Lemma~1.2.28 (b)]{DZ98}). Let us come to the details.

We proceed in two steps. First we consider the limit $ \delta\downarrow 0 $. Fix $ \eps >0 $ and a compact set $ \Lambda\subset\R^d $. Note that there is a compact set $ K_\Lambda\subset\Mcal_1(\Ccal) $ such that, for every $ \delta>0 $, the set $ F^{2\eps+\delta} $ can be replaced by $ F^{2\eps+\delta}\cap K_\Lambda $ without changing the value of the infimum on the left hand side of \eqref{countingupper8}. This can be seen as follows. From \eqref{exptightuniform}, together with the lower bound in the large deviations principle for $ L_N $ (see the proof of Proposition~\ref{proplowerpairmeasure}), one deduces that the set $ \{\mu\in\Mcal_1(\Ccal)\colon \inf_{\Ucal,\eta} I_{\Ucal}^{\ssup{\eta}}(\mu)\le C\} $ is compact for all $ C\in [0,\infty) $, where the infimum is taken over all partitions $ \Ucal $ of $ \Lambda $ and over all $ \eta\in\Mcal_1(\Sigma^2) $ (adapt the proof of \cite[Lemma~1.2.28 (b)]{DZ98}). Hence, also the set
$$
K_{\Lambda,C}:=\Bigl\{\mu\in\Mcal_1(\Ccal)\colon\inf_{\Ucal,\eta} \big\{H(\eta|m\otimes\eta^{\ssup{2}})+I_{\Ucal}^{\ssup{\eta}}(\mu)\big\}\le C\Bigr\}
$$ 
is compact. Choosing $ C $ large enough we can pick $ K_\Lambda=K_{\Lambda,C} $. 

For $\mu\in\Mcal_1(\Ccal)$ and $q\in\Mcal_1(\Lambda\times \Lambda)$, introduce 
\begin{equation}\label{JUcalqdef}
I_{\Ucal}^{\ssup{q}}(\mu)=\sup_{\Phi\in\Ccal_{\rm b}(\Ccal)}\Big[\langle \Phi,\mu\rangle-\int_\Lambda\int_\Lambda q(\d x,\d y)\,\log \E_{U_{r(x)},U_{r(x)}}^{\beta}\bigl[{\rm e}^{\Phi(B)}\bigr]\Big],
\end{equation} 
where $ r(x)\in\Sigma $ is defined by $ x\in U_{r(x)} $. Then we have,
\begin{equation}\label{countingupper9}
\begin{aligned}
\inf_{\mu\in  F^{2\eps+\delta}\cap K_\Lambda}&\inf_{\eta\in\Mcal_1^{\ssup{\eps}}(\Sigma^2)} \Bigl\{H(\eta|m_{\Ucal}\otimes\eta^{\ssup{2}})  +I_{\Ucal}^{\ssup{\eta}}(\mu)\Bigr\}\\
&\ge\inf_{\mu\in F^{2\eps+\delta}\cap K_\Lambda} \inf_{q\in\Mcal_1^{\ssup{\eps}}(\Lambda^2)}\Bigl\{H(q|\m_{\Ucal}\otimes q^{\ssup{2}})+I_{\Ucal}^{\ssup{q}}(\mu)\Big\},
\end{aligned}
\end{equation} 
as is seen from considering 
$$
q(\d x,\d y)=\sum_{r,s\in\Sigma} \frac{\eta(r,s)}{|U_r\times U_s|}\1_{U_r\times U_s}(x,y)\,\d x\d y.
$$
Here $ \m_{\Ucal}(\d x)=\sum_{r\in\Sigma}\frac{m(r)}{|U_r|}\1_{U_r}(x)\d x\in\Mcal_1(\Lambda) $.

Fix $ \Phi\in\Ccal_{\rm b}(\Ccal) $ and note that there is $C_{\Phi,\delta}>0$ satisfying $ \lim_{\delta\downarrow 0}C_{\Phi,\delta}=0 $ such that, for any $ q\in\Mcal_1(\Lambda^2) $,
\begin{equation}\label{Cdeltadef}
\begin{aligned}
\int_\Lambda\int_\Lambda q(\d x,\d y)\,\log \E_{U_{r(x)},U_{r(x)}}^{\beta}\bigl[{\rm e}^{\Phi(B)}\bigr]\le \int_\Lambda\int_\Lambda q(\d x,\d y)\,\log \E_{x,y}^{\beta}\bigl[{\rm e}^{\Phi(B)}\bigr] +C_{\Phi,\delta},
\end{aligned}
\end{equation} 
as follows from the proof of Lemma \ref{estimationsetendpoints} (we used that the fineness $ f_{\Ucal}$ is not larger than $ \delta$). 
We recall the representation of the entropy as a Legendre transform (see \cite[Lemma~3.2.13]{DS01}): for any $q\in\Mcal_1(\Lambda\times \Lambda)$,
\begin{equation}\label{HLegTrafo}
H(q|\m_{\Ucal}\otimes q^{\ssup{2}})=\sup_{g\in\Ccal_{\rm b}(\Lambda\times \Lambda)}\Big[\langle g,q\rangle-\log\langle {\rm e}^g,\m_{\Ucal}\otimes q^{\ssup{2}}\rangle\Big].
\end{equation}

We now write $\Ucal_\delta$  instead of $\Ucal$. For $\delta>0$, let $\mu_{\delta}\in F^{2\eps+\delta}\cap K_\Lambda $ and $ q_{\delta}\in \Mcal_1^{\ssup{\eps}}(\Lambda^2)$ be minimisers for the variational formula on the right hand side of \eqref{countingupper9}. Since $\Lambda$ is compact, as $ \delta\downarrow 0$, along suitable subsequences,  $\mu_{\delta}$ and  $q_{\delta}$  converge weakly towards suitable $\mu\in F^{1}\cap K_\Lambda $ and $q\in\Mcal_1(\Lambda^2)$,  respectively. In particular, $q^{\ssup{2}}_{\delta}$ converges weakly towards $q^{\ssup{2}}$ and $q^{\ssup{1}}_{\delta}$ converges weakly towards $q^{\ssup{1}}$. Certainly, we have $\mu\in \cap_{\delta>0}F^{2\eps+\delta}\cap K_\Lambda=F^{2\eps}\cap K_\Lambda $ since $F^{2\eps}$ is closed, and $q\in\Mcal_1^{\ssup{\eps}}(\Lambda^2)$. From \eqref{JUcalqdef}, \eqref{Cdeltadef} and \eqref{HLegTrafo}, we have, for any $ \delta>0$,
\begin{equation}\label{countingupper10}
\begin{aligned}
\inf_{\mu\in F^{2\eps+\delta}\cap K_\Lambda}& \inf_{q\in\Mcal_1^{\ssup{\eps}}(\Lambda^2)}\Bigl\{H(q|\m_{\Ucal_\delta}\otimes q^{\ssup{2}})+I_{\Ucal}^{\ssup{q}}(\mu)\Big\}\\
&\ge\langle g,q_{\delta}\rangle-\log\langle {\rm e}^g,\m_{\Ucal_\delta}\otimes q_{\delta}^{\ssup{2}}\rangle+\langle\Phi,\mu_{\delta}\rangle -\int_\Lambda\int _\Lambda q_{\delta}(\d x,\d y)\log\E_{x,y}^{\beta}\Bigl({\rm e}^{\Phi(B)}\Bigr)-C_{\Phi,\delta},
\end{aligned}
\end{equation} 
where $g\in\Ccal_{\rm b}(\Lambda\times \Lambda)$ and $\Phi\in\Ccal_{\rm b}(\Ccal)$ are arbitrary. Note that $\m_{\Ucal_\delta}\to \m_\Lambda$ weakly as $\delta\downarrow 0$, where $\m_\Lambda$ is the conditional distribution of $\m$ given $\Lambda$. Hence, $\m_{\Ucal_\delta}\otimes q_{\delta}^{\ssup{2}}$ converges, as $\delta\downarrow 0$, weakly towards $\m_\Lambda\otimes q^{\ssup{2}}$. Consequently, letting $ \delta\downarrow 0$ on the right hand side of \eqref{countingupper10} and recalling \eqref{countingupper9}, we obtain that 
\begin{equation}\label{countingupper12}
\begin{aligned}
\liminf_{\delta\downarrow 0}\inf_{\mu\in  F^{2\eps+\delta}\cap K_\Lambda}&\inf_{\eta\in\Mcal_1^{\ssup{\eps}}(\Sigma^2)} \Bigl\{H(\eta|m\otimes\eta^{\ssup{2}})  +I_{\Ucal}^{\ssup{\eta}}(\mu)\Bigr\}\\
&\ge \langle g,q\rangle-\log\langle {\rm e}^g,\m_{\Lambda}\otimes q^{\ssup{2}}\rangle+\langle\Phi,\mu\rangle -\int_\Lambda\int _\Lambda q (\d x,\d y)\, \log\E_{x,y}^{\beta}\bigl[{\rm e}^{\Phi(B)}\bigr].
\end{aligned}
\end{equation}
Since this holds for any $g\in\Ccal_{\rm b}(\Lambda\times \Lambda)$ and $\Phi\in\Ccal_{\rm b}(\Ccal)$, the left hand side is not smaller than $H(q|\m_\Lambda\otimes q^{\ssup{2}})+I^{\ssup{q}}(\mu)$, where we extended $q$ trivially to a probability measure on $\R^d\times \R^d$ (with support in $\Lambda\times \Lambda$). Hence, 
\begin{equation}\label{countingupper11}
\begin{aligned}
\mbox{l.h.s.~\eqref{countingupper12}}&\geq \inf_{\mu\in F^{2\eps}}\inf_{q\in \Mcal_1^{\ssup{\eps}}(\R^d\times \R^d)}\Big\{H(q|\m_\Lambda\otimes q^{\ssup{2}})+I^{\ssup q}(\mu)\Big\}.
\end{aligned}
\end{equation}

In the second step of the proof, we replace $ \eps $ by $ \eps_N\downarrow 0 $ and $ \Lambda $ by $ \Lambda_N\uparrow\R^d $ and consider the limit as $ N\to\infty $. Clearly, $\m_{\Lambda_N}\to \m$ weakly. For any $N\in\N$, pick $\mu_{N}\in F^{2\eps_N} $ and $q_N\in \Mcal_1^{\ssup{\eps_N}}(\R^d\times \R^d)$ such that the sequence $ (H(q_N|\m_{\Lambda_N}\otimes q^{\ssup{2}}_N)+I^{\ssup{q_N}}(\mu_N))_{N\in\N} $ converges to the left hand side of \eqref{countingupper8} and may therefore be assumed to be bounded. Since
$$
H(q_N|\m_{\Lambda_N}\otimes q_N^{\ssup{2}})=H(q^{\ssup{1}}_N|\m_{\Lambda_N})+H(q_N|q_N^{\ssup{1}}\otimes q_N^{\ssup{2}}),
$$
the sequence $ (H(q^{\ssup{1}}_N|\m_{\Lambda_N}))_{N\in\N} $ is also bounded.  Since $ H(q^{\ssup{1}}_N|\m)=H(q^{\ssup{1}}_N|\m_{\Lambda_N}) $,  the sequence $ (q^{\ssup{1}}_N)_{N\in\N} $ is tight, because the level sets of the relative entropy are compact (see \cite[Lemma~6.2.12]{DZ98}). 
Since $ \d(q^{\ssup{1}}_N,q^{\ssup{2}}_N)\le 2\eps_N\to 0 $ as $ N\to\infty $, also $ (q^{\ssup{2}}_N)_{N\in\N} $ is tight. By boundedness of $(H(q_N|q_N^{\ssup{1}}\otimes q_N^{\ssup{2}}))_{N\in\N}$, also the set $ Q:=\{q_N\colon N\in\N\} $ is tight. According to Prohorov's theorem, we may assume that $ q_N\Rightarrow q  $ for some $ q\in\Mcal_1(\R^d\times\R^d)$. Since also $ q^{\ssup{1}}_N\Rightarrow q^{\ssup{1}} $ and $  q^{\ssup{2}}_N\Rightarrow q^{\ssup{2}} $ and $ \d(q^{\ssup{1}}_N,q^{\ssup{2}}_N)\to 0 $, we have that $  q\in\Mcal_1^{\ssup{\rm s}}(\R^d\times\R^d)$. 

For sufficiently large $ C>0 $, the sequence  $ (\mu_N)_{N\in\N} $ is contained in the set $ \{\mu\in\Mcal_1(\Ccal)\colon\inf_{N\in\N}I^{\ssup{q_N}}(\mu)\le C\} $. It turns out that this set is relatively compact. For proving this, it suffices to find a family of compact sets $ K_L\subset\Ccal$, $L>0 $, such that
$$
\lim_{L\to\infty}\inf_{q\in Q}\inf_{K_L^{\rm c}} I^{\ssup{q}}=\infty.
$$
Consider a sequence of compact sets $ \Lambda_N\uparrow\R^d $ (not necessarily those we picked above) and a sequence of partitions $ \Ucal_N =\{U_r\colon r\in\Sigma_N\}$ of $ \Lambda_N$ whose fineness vanishes as $ N\to\infty $. Given $ q\in Q $, pick $\eta_N^{\ssup q}\in\Mcal_1^{\ssup{N}}(\Sigma_N^2) $ such that the probability measures
$$
q_N^{\ssup q}(\d x,\d y):=\sum_{r,s\in\Sigma_N}\frac{\eta_N^{\ssup q}(r,s)}{|U_r\times U_s|}\1_{U_r\times U_s}(x,y)\,\d x\d y
$$
converge weakly to $ q $. Then the sequence of empirical path measures, $ (L_N)_{N\in\N} $, is exponentially tight under $ P^\beta_{\eta_N^{\ssup q},N,\Ucal_N}$, uniformly in $ q\in Q $ (see Lemma~\ref{tightnessPeta}). Furthermore, it satisfies a large deviations principle with rate function $ I^{\ssup{q}} $. This  is seen as follows. The logarithmic moment generating function of $L_N$ under 
$ P^\beta_{\eta_N^{\ssup q},N,\Ucal_N}$, defined in \eqref{logmomgenfctLN}, is easily shown to converge towards the function $\Phi\mapsto\int_{\R^d}\int_{\R^d}q(\d x,\d y)\log \E^\beta_{x,y}[{\rm e}^{\Phi(B)}]$. Since its Fenchel-Legendre transform is equal to $I^{\ssup q}$, the G\"artner-Ellis theorem implies the mentioned large deviations principle.
 
For $L\in\N$, pick a compact set $ K_L\subset\Ccal $ such that $ P^\beta_{\eta_N^{\ssup q},N,\Ucal_N}(L_N\in K_L^{\rm c})\le{\rm e}^{-NL} $ for all $ L,N\in\N $ and $q\in Q$. Using the lower bound in the mentioned large deviations principle, this implies that 
$$
\inf_ {q\in Q}\inf_{K_L^{\rm c}} I^{\ssup{q}}\ge -\liminf_{N\to\infty}\frac{1}{N}\log P^\beta_{\eta_N,N,\Ucal_N}(L_N\in K_L^{\rm c})\ge L.
$$
Hence, the sequence $ (\mu_N)_{N\in\N} $ is tight.

Therefore, we may assume that $ \mu_N\Rightarrow \mu $ for some $ \mu\in F^{\ssup{1}}$. Since $\mu_N\in F^{\ssup{2\eps_N}}$ for any $N\in\N$ and since $\eps_N\to 0$, we even have that $\mu\in F$, since $F$ is closed. Now in the same way as we derived \eqref{countingupper11}, one derives that \eqref{countingupper8} holds.
\qed
\end{proofsect}


\subsection{Exponential tightness.}\label{finishproofmain}

\noindent In this section, we prove the necessary exponential tightness assertions for the sequence of the empirical path measures under the symmetrised measures,  $ P_{\m,N}^{\ssup{\rm sym}} $, and  under the mixed product measures, $  P_{\eta_N,N,\Ucal_N}^{\beta} $. The proof of the latter exponential tightness is a variant of the standard proof for laws of empirical measures. Here, the main ingredient is the product structure of the probability measure. The proof of the first exponential tightness exploits a compactification argument due to the starting distribution $ \m\in\Mcal_1(\R^d) $.

\begin{lemma}\label{exptight}
Let $ \m\in\Mcal_1(\R^d) $ be the initial distribution. Then the family of distributions of the empirical path measures $ L_N $ under the symmetrised measure $ \P_{\m,N}^{\ssup{\rm sym}} $ is exponentially tight.
\end{lemma}
\begin{proofsect}{Proof}
The proof is in the spirit of the proof of \cite[Lemma~6.2.6]{DZ98}. 
For $ l\in\N $, choose a box $ Q_l\subset\R^d $ such that $ \m(Q_l^{\rm c})\le {\rm e}^{-l^2}. $ Furthermore, choose $ \delta_l > 0 $ so small that
\begin{equation}
\sup_{x,y\in Q_l}\P_{x,y}^{\beta}\Bigl(\sup_{|s-t|\le \delta_l}|B_s-B_t| > \frac{1}{l}\Bigr)\le {\rm e}^{-l^2}.
\end{equation}
Consider  
$$ 
A_l=\Big\{f\in\Ccal\colon f(0)\in Q_l, f(\beta)\in Q_l, \sup_{|s-t|\le \delta_l}|f(s)-f(t)|\le \frac{1}{l}\Big\}. 
$$ 
According to Arzel\`a-Ascoli's theorem, $ A_l $ is relative compact in $ \Ccal$. Put $ M_l:=\{\mu\in\Mcal_1(\Ccal)\colon\mu(\overline{A}_l^{\rm c})\le \frac{1}{l}\} $ and note that $ M_l $ is closed by Portmanteau's theorem. Let $ L\in\N $ be given and consider $ K_L:=\bigcap_{l=L}^{\infty} M_l $. It is easy to see that $K_L$ is tight, hence $\overline K_L$ is compact by Prohorov's theorem. We shall show that $ \P_{\m,N}^{\ssup{\rm sym}}(L_N\in \overline K_L^{\rm c})\le {\rm e}^{-LN} $ for any $N\in\N$, which implies the assertion. Observe that 
$$
\begin{aligned}
\{L_N\in M^{\rm c}_l\}&\subset\Big\{\,\sharp\{i\in\{1,\ldots,N\}\colon B^{\ssup{i}}\in A^{\rm c}_l\} >\frac{N}{l}\Big\}\\&\subset \Big\{\sharp\{i\colon B^{\ssup{i}}_0\in Q_l^{\rm c}\}\ge \frac{N}{3l}\Big\}\cup\Big\{\sharp\{i\colon B^{\ssup{i}}_{\beta}\in Q_l^{\rm c}\}\ge\frac{N}{3l}\Big\}\\&\quad\cup  \Big\{\sharp\Big\{i\colon B^{\ssup{i}}_0\in Q_l, B^{\ssup{i}}_{\beta}\in Q_l, \sup_{|s-t|\le \delta_l}|B^{\ssup{i}}_s-B^{\ssup{i}}_t| >\frac{1}{l}\Big\}\ge \frac{N}{3l}\Big\}.
\end{aligned}
$$
Clearly,
\begin{equation}
\begin{aligned}
\P_{\m,N}^{\ssup{\rm sym}}\Big(\sharp\{i\colon B^{\ssup{i}}_{\beta}\in Q_l^{\rm c}\}\ge \frac{N}{3l}\Big)&=\P_{\m,N}^{\ssup{\rm sym}}\Big(\sharp\{i\colon B^{\ssup{i}}_{0}\in Q_l^{\rm c}\}\ge \frac{N}{3l}\Big)\\
& \le \sum_{\heap{I\subset\{1,\ldots,N\}\colon}{|I|\ge \frac{N}{3l}}}\frac{1}{N!}\sum_{\s\in\Sym_N}\int_{(\R^d)^N}\prod_{i=1}^N\m(\d x_i)\bigotimes_{i=1}^N\P_{x_i,x_{\s(i)}}\big(\forall\; i\in I\colon B^{\ssup{i}}_0\in Q^{\rm c}_l\big)\\
&\le\sum_{|I|\ge \frac{N}{3l}} \m(Q_l^{\rm c})^{|I|}\le 2^N{\rm e}^{-lN/3}.
\end{aligned}
\end{equation}
Furthermore,
\begin{equation}
\begin{aligned}
\P&_{\m,N}^{\ssup{\rm sym}}\Bigl(\sharp\Bigl\{i\colon B^{\ssup{i}}_0\in Q_l,B^{\ssup{i}}_{\beta}\in Q_l, \sup_{|s-t|\le\delta_l}|B^{\ssup{i}}_0-B^{\ssup{i}}_{\beta}| >\frac{1}{l}\Bigr\}\ge \frac{N}{3l}\Bigr)\\
& \le \sum_{\heap{I\subset\{1,\ldots,N\}\colon}{|I|\ge \frac{N}{3l}}}\frac{1}{N!}\sum_{\s\in\Sym_N}\int_{(\R^d)^N}\prod_{i=1}^N\m(\d x_i)\\
&\qquad\qquad\bigotimes_{i=1}^N\P_{x_i,x_{\s(i)}}^{\beta}\Bigl(\forall\; i\in I\colon B^{\ssup{i}}_0\in Q_l,B^{\ssup{i}}_{\beta}\in Q_l,\sup_{|s-t|\le\delta_l}|B^{\ssup{i}}_s-B_t^{\ssup{i}}| > \frac{1}{l}\Bigr)\\
&\le \sum_{|I|\ge \frac{N}{3l}}\sup_{(y_i)_{i\in I}\in Q_l^{I}}\int_{Q_l^{I}}\prod_{i\in I}\m(\d x_i)\prod\limits_{i\in I}\P_{x_i,y_i}^{\beta}\Bigl(\sup_{|s-t|\le\delta_l}|B^{\ssup{i}}_s-B_t^{\ssup{i}}| > \frac{1}{l}\Bigr)\\
&\le \sum_{|I|\ge \frac{N}{3l}}{\rm e}^{-lN/3}\le 2^N{\rm e}^{-lN/3}.
\end{aligned}
\end{equation} Hence, 
$$ 
\P_{\m,N}^{\ssup{\rm sym}}(L_N\in K_L^{\rm c})\le\sum_{l=L}^{\infty}\P_{\m,N}^{\ssup{\rm sym}}(L_N\in M^{\rm c}_l)\le 3\times2^N\sum_{l=L}^{\infty}{\rm e}^{-lN/3}\le 6\times2^N{\rm e}^{-NL/3}\le {\rm e}^{-NL/5}
$$ for all large $ N $ if $ L> 24 $. This ends the proof.
\qed\end{proofsect} 

Now we prove the exponential tightness of the empirical path measures $L_N$ under the measures $P_{\eta,N,\Ucal}^{\beta}$ introduced in \eqref{low3}. We continue to use the notation introduced at the beginning of Section~\ref{sec-lowproof}.

\begin{lemma}\label{tightnessPeta}
Let $ (\Lambda_N)_{N\in\N} $ be a sequence of compact subsets of $ \R^d $ and let $ (\Ucal_N)_{N\in\N} $ be a sequence of partitions $ \Ucal_N=\{U_r\colon r\in\Sigma_N\} $ of $ \Lambda_N $. For any $ N\in\N $, let $ \eta_N$ be in $\Mcal_1^{\ssup{N}}(\Sigma_N^2) $ such that the sequence of probability measures $ q_N $ defined by
$$
q_N(\d x,\d y)=\sum_{r,s\in\Sigma_N}\frac{\eta_N(r,s)}{|U_r\times U_s|}\1_{U_r\times U_s}(x,y)\,\d x\d y, 
$$ 
is tight. Then the families of distributions of the empirical path measures $ L_N $ and the one of the means $ Y_N $ of occupation measures under the measures $ P_{\eta_N,N,\Ucal_N}^{\beta} $ are exponentially tight.
\end{lemma}

\begin{proofsect}{Proof}
We prove the exponential tightness for the empirical path measures, the one for the means $ Y_N $ follows analogously. As we have seen at the beginning of the proof of Lemma~\ref{exptight}, for any $ l\in\N $ there exists a compact set $ Q_l\subset\R^d $ such that for all $ N\in\N $ we have $ q_N((Q_l\times Q_l)^{\rm c})\le \frac{1}{6l} $.  Furthermore, there exists a compact set $ \Gamma_l\subset\Ccal $ such that
\begin{equation}
\sup_{x,y\in Q_{l+1}}\P_{x,y}^{\beta}(B\in\Gamma_l^{\rm c})\le {\rm e}^{-2l^2}({\rm e}^l-1).
\end{equation}
The set $ M_l=\{\nu\in\Mcal_1(\Ccal)\colon \nu(\Gamma_l^{\rm c})\le 1/l\} $ is closed by  Portmanteau's theorem. For $ L\in\N $ define $ K_L=\bigcap_{l=L}^{\infty}M_l $. By Prohorov's theorem, each $ K_L $ is a relative compact subset of $ \Mcal_1(\Ccal) $. We may assume that $ \diam\, U_r <1 $ for any $ r\in\Sigma_N $. Then Chebycheff's inequality gives that for any $ N\in\N $, any partition $ \Ucal_N $ of $ \Lambda_N $ and any $ \eta_N\in\Mcal_1^{\ssup{N}}(\Sigma^2_N) $  
\begin{equation}
\begin{aligned}
P_{\eta_N,N,\Ucal_N}^{\beta}(L_N\notin M_l)&=P_{\eta_N,N,\Ucal_N}^{\beta}\Bigl(L_N(\Gamma_l^{\rm c})>\frac{1}{l}\Bigr)\le E_{\eta_N,N,\Ucal_N}^{\beta}\Bigl[{\rm e}^{2Nl^2(L_N(\Gamma_l^{\rm c})-1/l)}\Bigr]\\
&= {\rm e}^{-2Nl}E_{\eta_N,N,\Ucal_N}^{\beta}\Bigl[\exp\Bigl(2l^2\sum_{i=1}^N\1\{B^{\ssup{i}}\in \Gamma_l^{\rm c}\}\Bigr)\Bigr]\\
&={\rm e}^{-2Nl}\prod_{r,s\in\Sigma_N}\E_{U_r,U_s}^{\beta}\Bigr[\exp\Bigr(2l^2\1\{B\in \Gamma_l^{\rm c}\}\Bigr)\Bigr]^{N\eta_N(r,s)}\\
&={\rm e}^{-2Nl}\prod_{r,s\in\Sigma_N}\Bigl(\P_{U_r,U_s}^{\beta}(B\in\Gamma_l)+{\rm e}^{2l^2}\P_{U_r,U_s}^{\beta}(B\in\Gamma_l^{\rm c})\Bigr)^{N\eta_N(r,s)}\\
& \le {\rm e}^{-2Nl}\Big(\prod_{\heap{r,s\in\Sigma_N\colon}{U_r\times U_s\subset Q_{l+1}^2}}
({\rm e}^{l})^{N\eta_N(r,s)}\Big)({\rm e}^{3l^2})^{Nq_N((Q_l\times Q_l)^{\rm c})}\le {\rm e}^{-Nl/2},
\end{aligned}
\end{equation}
where in the last line we also used that (see \eqref{PUUdef}) 
$$
\P_{U_r,U_s}^{\beta}(B\in\Gamma_l^{\rm c})\le{\rm e}^{-2l^2}({\rm e}^l-1)\quad\mbox{ if }\,U_r\times U_s\subset Q_{l+1}^{2},
$$
that $ 1+{\rm e}^{2l^2}\le {\rm e}^{3l^2} $ and that $q_N((Q_l\times Q_l)^{\rm c})\le \frac{1}{6l} $.
Therefore,
\begin{equation}
P_{\eta_N,N,\Ucal_N}^{\beta}(L_N\notin K_L)\le\sum_{l=L}^{\infty}P_{\eta_N,N,\Ucal_N}^{\beta}(L_N\notin M_l)\le \sum_{l=L}^{\infty}{\rm e}^{-Nl/2}\le 2{\rm e}^{-NL/2},
\end{equation}
which implies the exponential tightness.
\end{proofsect}
\qed

\subsection[Proof of Theorem~\ref{main2}]{Proof of Theorem~\ref{main2}}\label{sec-proofYN}

\noindent In this section we prove Theorem \ref{main2}. We recall that a large-deviation principle for $Y_N$ under $\P^{\ssup{\rm sym}}_{\m,N}$ with rate function $ \widetilde J^{\ssup{\rm sym}}_{\m} $ (see \eqref{ratefunctioncontract}) directly follows from the principle of Theorem~\ref{main} for $L_N$ via the contraction principle \cite[Th.~4.2.1]{DZ98}, since $Y_N=\Psi(L_N) $, where $\Psi(\mu)=\frac 1\beta\int_0^\beta \mu\circ \pi_s^{-1}\,\d s$. The rate function is given as $ \widetilde J^{\ssup{\rm sym}}_{\m}$ defined in \eqref{ratefctJtilde}. Therefore, it suffices to show that $ \widetilde J^{\ssup{\rm sym}}_{\m}$ coincides with $J^{\ssup{\rm sym}}_{\m} $ introduced in \eqref{Jsymdef}. For this, it suffices to show that the two functions $\widetilde J^{\ssup{q}}$ and $J^{\ssup{q}}$, defined in \eqref{ratefunctioncontract} and \eqref{Jsymqdef}, coincide for any $ q\in\Mcal_1^{\ssup{\rm s}}(\R^d\times\R^d) $.

Fix $ q\in\Mcal_1^{\ssup{\rm s}}(\R^d\times\R^d) $ and let us first show that $ \widetilde J^{\ssup{q}}\ge J^{\ssup{q}}$. Given $ \mu\in\Mcal_1(\Ccal) $, we specialise the supremum over $ \Phi\in\Ccal_{\rm b}(\Ccal) $ in the definition \eqref{Jqdef} of $ I^{\ssup{q}} $ to functions of the form $ \Phi(\omega)=\frac{1}{\beta}\int_0^{\beta}\d s\, f(\omega(s)) $ with $ f\in\Ccal_{\rm b}(\R^d) $, to obtain that
\begin{equation}\label{JtildeqJq}
\begin{aligned}
I^{\ssup{q}}(\mu)&\ge \sup\limits_{f\in\Ccal_{\rm b}(\R^d)}\Big\{\int_{\Ccal}\mu(\d \omega)\frac{1}{\beta}\int_0^{\beta}\d s\,f(\omega(s))-\int_{\R^d}\int_{\R^d}q(\d x,\d y)\log\E_{x,y}^{\beta}\big[{\rm e}^{\frac{1}{\beta}\int_0^{\beta}  f(B_s)\,\d s}\big]\Big\}\\
&= J^{\ssup{q}}(\Psi(\mu)).
\end{aligned}
\end{equation}
Taking the infimum over all $\mu$ satisfying $\Psi(\mu)=p$, it is clear that $ \widetilde J^{\ssup{q}}(p)\ge J^{\ssup{q}}(p)$ for any $p\in\Mcal_1(\R^d)$.

It remains to show the complementary bound, $ \widetilde J^{\ssup{q}}(p)\le J^{\ssup{q}}(p)$. Proving this directly in an analytical way seems to cause major difficulties. Therefore, we proceed in an indirect way by showing that both $\widetilde J^{\ssup{q}}$ and $J^{\ssup{q}}$ are the rate function for the same large deviations principle. By the uniqueness of the rate function, this implies the assertion (even without using \eqref{JtildeqJq}).

Measures that satisfy a large deviations principle with rate function $I^{\ssup{q}}$ have been constructed at the end of the proof of Lemma~\ref{LDPargupper}. Indeed, consider a sequence of compact sets $ \Lambda_N\uparrow\R^d $ and a sequence of partitions $ \Ucal_N =\{U_r\colon r\in\Sigma_N\}$ of $ \Lambda_N$ whose fineness vanishes as $ N\to\infty $. Pick $\eta_N\in\Mcal_1^{\ssup{N}}(\Sigma_N^2) $ such that the probability measures
$$
q_N(\d x,\d y):=\sum_{r,s\in\Sigma_N}\frac{\eta_N(r,s)}{|U_r\times U_s|}\1_{U_r\times U_s}(x,y)\,\d x\d y
$$
converge weakly to $ q $. According to Lemma~\ref{tightnessPeta}, the sequence of empirical path measures, $ (L_N)_{N\in\N} $, is exponentially tight under $ P^\beta_{\eta_N,N,\Ucal_N}$. As has been explained in the proof of Lemma~\ref{LDPargupper}, it satisfies a large deviations principle with rate function $ I^{\ssup{q}} $. According to the contraction principle, the sequence $(Y_N)_{N\in\N} $ satisfies, under the measures $ P^\beta_{\eta_N,N,\Ucal_N}$, a large deviations principle with rate function $ \widetilde J^{\ssup{q}} $.

Now we show that $(Y_N)_{N\in\N} $
satisfies, under the measures $ P^\beta_{\eta_N,N,\Ucal_N}$, a large deviations principle with rate function $ J^{\ssup{q}} $, which ends the proof. For this, we have to consider the logarithmic moment generating function of $Y_N$ under $ P^\beta_{\eta_N,N,\Ucal_N}$, which is identified, for any $ f\in\Ccal_{\rm b}(\R^d) $, as
\begin{equation}
\begin{aligned}
\Lcal_N(f)&:=\log \E^{\beta}_{\eta_N,N,\Ucal_N}\big[{\rm e}^{N\langle f ,Y_N\rangle}\big]
=\log\Bigl(\prod_{r,s\in\Sigma_N}\E^{\beta}_{U_r,U_s}\big[{\rm e}^{\int_0^{\beta}f(B_s)\,\d s}\big]^{N\eta_N(r,s)}\Bigr)\\
&=N\sum_{r,s\in\Sigma_N}\eta_N(r,s)\log \E^{\beta}_{U_r,U_s}\big[{\rm e}^{\int_0^{\beta}f(B_s)\,\d s}\big]\\
&=N\int_{\R^d}\int_{\R^d}q_N(\d x,\d y)\,\log\E^{\beta}_{U_{r_N(x)},U_{r_N(y)}}\big[{\rm e}^{\int_0^{\beta}f(B_s)\,\d s}\big],
\end{aligned}
\end{equation} 
where $ r_N(x)\in\Sigma_N $ is defined by $ x\in U_{r_N(x)}$. From the proof of Lemma~\ref{estimationsetendpoints} it is seen that 
$$ 
\lim_{N\to\infty} \E^{\beta}_{U_{r_N(x)},U_{r_N(y)}}\big[{\rm e}^{\int_0^{\beta}f(B_s)\,\d s}\big]=\E^{\beta}_{x,y}\big[{\rm e}^{\int_0^{\beta}f(B_s)\,\d s}\big], 
$$ 
uniformly in $ x,y $ on compact sets. Recall that $ q_N\to q $ as $ N\to\infty $ weakly. Hence, the limit $ \Lcal(f)=\lim_{N\to\infty}\frac 1N\Lcal_N(f) $ exists, and 
$$
\Lcal(f)=\int_{\R^d}\int_{\R^d}q(\d x,\d y)\,\log\E^{\beta}_{x,y}\big[{\rm e}^ {\int_0^{\beta}f(B_s)\,\d s}\big].
$$ 
It is easily seen that $ \Lcal $ is lower semi continuous and  G\^ateaux differentiable. Note further that the Fenchel-Legendre transform of $\Lcal$ is equal to $J^{\ssup q}$. Furthermore, according to Lemma~\ref{tightnessPeta}, the sequence $(Y_N)_{N\in\N} $ is exponentially tight under  $ (P_{\eta_N,N,\Ucal_N}^{\beta})_{N\in\N} $. Hence, the G\"artner-Ellis theorem \cite[4.5.27]{DZ98} implies that $(Y_N)_{N\in\N} $
satisfies, under the measures $ P^\beta_{\eta_N,N,\Ucal_N}$, a large deviations principle with rate function $ J^{\ssup{q}} $, which ends the proof.

\section{Appendix: large deviations}\label{sec-LDP}

\noindent For the convenience of our reader, we repeat the notion of a large-deviation principle and of the most important facts that are used in the present paper. See \cite{DZ98} for a comprehensive treatment of this theory. 

Let $ \Xcal $ denote a topological vector space.  A lower semi-continuous function $ I\colon \Xcal\to [0,\infty] $ is called a {\it rate function\/} if  $ I $ is not identical $ \infty$ and  has compact level sets, i.e., if $ I^{-1}([0,c])=\{x\in\Xcal\colon I(x)\le c\} $ is compact for any $ c\ge 0 $. A sequence $(X_N)_{N\in\N}$ of $\Xcal$-valued random variables $X_N$  satisfies the {\it large-deviation upper bound\/} with {\it speed\/} $a_N$ and rate function $I$ if, for any closed subset $F$ of $\Xcal$,
\begin{equation}\label{LDPupper}
\limsup_{N\to\infty}\frac 1{a_N}\log \P(X_N\in F)\leq -\inf_{x\in F}I(x),
\end{equation}
and it satisfies the {\it large-deviation lower bound\/} if, for any open subset $G$ of $\Xcal$,
\begin{equation}\label{LDPlower}
\liminf_{N\to\infty}\frac 1{a_N}\log \P(X_N\in G)\leq -\inf_{x\in G}I(x).
\end{equation}
If both, upper and lower bound, are satisfied, one says that $(X_N)_N$ satisfies a {\it large-deviation principle}. The principle is called {\it weak\/} if the upper bound in \eqref{LDPupper} holds only for {\it compact\/} sets $F$. A weak principle can be strengthened to a full one by showing that the sequence of distributions of $X_N$ is {\it exponentially tight}, i.e., if for any $L>0$ there is a compact subset $K_L$ of $\Xcal$ such that $\P(X_N\in K_L^{\rm c})\leq {\rm e}^{-LN}$ for any $N\in\N$. 

One of the most important conclusions from a large deviation principle is {\it Varadhan's Lemma}, which says that, for any bounded and continuous function $F\colon \Xcal\to\R$,
$$
\lim_{N\to\infty}\frac 1N\log \int {\rm e}^{N F(X_N)}\,\d\P=-\inf_{x\in \Xcal}\big(I(x)-F(x)\big).
$$

All the above is usually stated for probability measures $\P$ only, but the notion easily extends to {\it sub}-probability measures $\P=\P_N$ depending on $N$. Indeed, first observe that the situation is not changed if $\P$ depends on $N$, since a large deviation principle depends only on distributions. Furthermore, the connection between probability distributions $\widetilde \P_N$ and sub-probability measures $\P_N$ is provided by the transformed measure $\widetilde \P_N(X_N\in A)=\P_N(X_N\in A)/\P_N(X_N\in\Xcal)$: if the measures $\P_N\circ X_N^{-1}$ satisfy a large deviation principle with rate function $I$, then the probability measures $\widetilde \P_N\circ X_N^{-1}$ satisfy a large deviation principle with rate function $I-\inf I$.

One standard situation in which a large deviation principle holds is the case where $\P$ is a probability measure, and $X_N=\frac 1N(Y_1+\dots+Y_N)$ is the mean of $N$ i.i.d.~$\Xcal$-valued random variables $Y_i$ whose moment generating function $M(F)=\int {\rm e}^{F(Y_1)}\,\d\P$ is finite for all elements $F$ of the topological dual space $\Xcal^*$ of $\Xcal$. In this case, the abstract Cram\'er theorem provides a weak large deviation principle for $(X_N)_{N\in\N}$ with rate function equal to the Legendre-Fenchel transform of $\log M$, i.e., $I(x)=\sup_{F\in \Xcal^*}(F(x)-\log M(F))$. An extension to independent, but not necessarily identically distributed random variables is provided by the abstract G\"artner-Ellis theorem.

In our large deviations results we shall rely on the following conventions. For $X=\Ccal$ or $X=\R^d$, we conceive $ \Mcal_1(X)$  as a closed convex subset of the
space $\Xcal=\Mcal(X)$ of all finite signed Borel measures on $X$. This is a
topological Hausdorff vector space whose topology is induced by the set
$\Ccal_{\rm b}(X)$ of all continuous bounded functions  $X\to\R$. Then $ \Ccal_{\rm b}(X)$ is  the topological dual of $ \Mcal(X)$ \cite[Lemma~3.2.3]{DS01}. When we speak of a
large deviation principle for $\Mcal_1(X)$-valued random variables, then we
mean a principle on $\Mcal(X)$ with a rate function that is tacitly
extended from $\Mcal_1(X)$ to $\Mcal(X)$ with the value $+\infty$.

\end{document}